\documentclass[12pt,a4paper]{article}
\usepackage{fullpage}
\usepackage{mdwlist}
\usepackage{amsmath}
\usepackage{amsthm}
\usepackage{amssymb}
\usepackage{mathrsfs}
\usepackage{amsfonts}
\usepackage{color}
\usepackage{tikz}
\usepackage[all]{xy}
\usepackage{enumitem}
\newcommand{\nd}{\triangleleft}
\newcommand{\R}{\ensuremath{\mathbb{R}}}
\newcommand{\Q}{\ensuremath{\mathbb{Q}}}
\newcommand{\Z}{\ensuremath{\mathbb{Z}}}
\newcommand{\N}{\ensuremath{\mathbb{N}}}

\newcommand{\C}{\ensuremath{\mathbb{C}}}
\newcommand{\F}{\ensuremath{\mathbb{F}}}

\newcommand{\co}{\ensuremath{\mathcal{O}}}
\newcommand{\cn}{\ensuremath{\mathcal{N}}}

\newcommand{\gl}{\ensuremath{\mathfrak{g}}}
\newcommand{\Gl}{\ensuremath{\mathfrak{gl}}}
\newcommand{\zl}{\ensuremath{\mathfrak{z}}}
\newcommand{\nl}{\ensuremath{\mathfrak{n}}}
\newcommand{\onl}{\ensuremath{\overline{\mathfrak{n}}}}
\newcommand{\ml}{\ensuremath{\mathfrak{m}}}
\newcommand{\tl}{\ensuremath{\mathfrak{t}}}
\newcommand{\pl}{\ensuremath{\mathfrak{p}}}
\newcommand{\bl}{\ensuremath{\mathfrak{b}}}

\newcommand{\coker}{\mathrm{coker}\,}
\renewcommand{\hom}{\mathrm{Hom}}
\newcommand{\eps}{\epsilon}
\newcommand{\weglaten}[1]{}
\newtheorem{thm}{}
\newtheorem{ste}[thm]{Theorem}
\newtheorem{mdef}[thm]{Definition}
\newtheorem{lem}[thm]{Lemma}
\newtheorem{gev}[thm]{Corollary}
\newtheorem{ver}[thm]{Conjecture}
\newtheorem{pro}[thm]{Proposition}

\begin{document}
\begin{center}
{\large\bf Nilpotent orbits:\\ finiteness, separability and Howe's conjecture}\\
\vspace{0.3cm}
{\large Julius Witte}\\
\vspace{0.1cm}
Radboud Universiteit Nijmegen\\
Heyendaalseweg 135, 6525AJ Nijmegen, the Netherlands\\
email: J.Witte@math.ru.nl
\end{center}

\section*{Abstract}
This paper is about nilpotent orbits of reductive groups over local non-Archimedean fields. In this paper we will try to identify for which groups there are only finitely many nilpotent orbits, for which groups the nilpotent orbits are separable and for which groups Howe's conjecture holds. For general reductive groups we get some partial results. For split reductive groups we get a classification in terms of the root data and the characteristic of the underlying local field. 
\tableofcontents
\vspace{0.3cm}
\noindent{\bf Acknowledgements}\\
The author would like to thank M. Solleveld for many useful discussions. The author is financially supported by NWO grant nr 613.009.033 and EW wiskundecluster GQT positie.

\section{Introduction}
Let $\F$ be a local non-Archimedean field of characteristic $p$ and $\mathcal{G}$ a connected reductive group defined over $\F$. In this paper we investigate the relation between the following statements for $G=\mathcal{G}(\F)$ (we will clarify the first six statements in \S2):
\begin{enumerate*}
 \item $p$ is good
 \item $p$ is very good
 \item $p$ does not divide the virtual number of components of $Z(G)$
 \item $p$ does not divide the virtual order of $\pi_1(G_{\text{der}})$
 \item all the nilpotent orbits are separable
 \item the regular nilpotent orbit is separable
 \item the number of nilpotent orbits is finite
 \item Howe's conjecture holds for $G$
\end{enumerate*}
If char $\F=0$ (including $\F=\C,\R$), then all these statements hold for $G$. In case $\F$ has positive characteristic these statements depend on $G$ and $p$.\\

For general $G$ we will prove the following implications
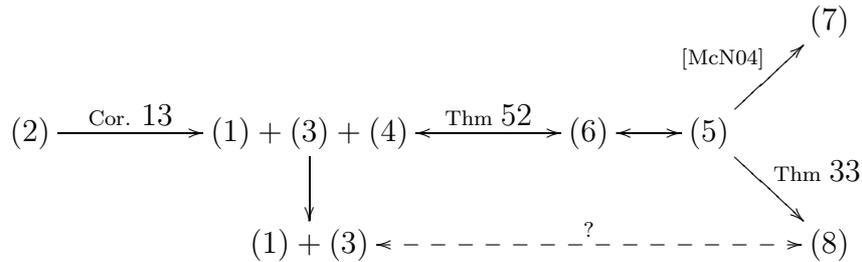
\[
 \xymatrix{&&&&&&(7) \\(2) \ar@{->}[rr]^-{\text{Cor. }\ref{2implies134}} && (1)+(3)+(4)\ar[rr]^-{\text{Thm } \ref{cnos}}\ar[d] && (6) \ar[ll]\ar[r]& (5) \ar[l]\ar[rd]^-{\text{Thm } \ref{Higc}}\ar[ru]^-{\text{\cite{Mc04}}}\\&&(1)+(3)\ar@{<-->}[rrrr]^-{?}&&&&(8)
 }
\]
The question if $(7)$ and $(5)$ are equivalent and the question if $(1)+(3)$ is equivalent to $(8)$ are still open.\\
If moreover $G$ is $\F$-split, then we get the following implications
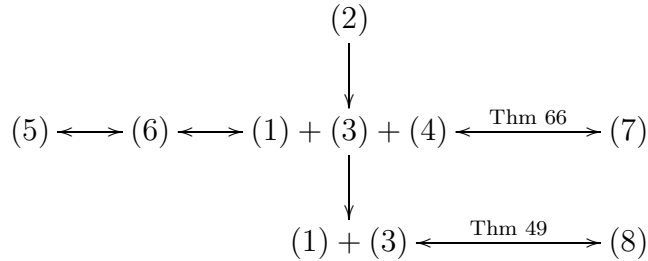
\[
 \xymatrix{ &&(2)\ar[d]&&\\(5)\ar[r] & (6) \ar[l]\ar[r]& (1)+(3)+(4) \ar[l]\ar[rr]^-{\text{Thm \ref{cfnisc}}}\ar[d] && (7)\ar[ll]\\ &&(1)+(3)\ar[rr]^-{\text{Thm \ref{Hscc}}}&&(8)\ar[ll]& }
\]
Besides the proofs of these implications we will also give counter examples for the non-implications. That $(1)+(3)$ does not imply $(4)$ can be seen by the example $PGL_p$. That $(1)+(3)+(4)$ does not imply $(2)$ can be seen by the example $GL_p$. That $(1)+(4)$ does not imply $(3)$ can be seen by the example $SL_p$. That $(3)+(4)$ does not imply $(1)$ can be seen in the simple groups of exceptional type. Thus for $\F$-split groups we have determined all the implications and non-implications between every possible combination of these 8 properties.\\

The first $4$ statements are related to $p$ and the root datum of $G$ and the last $4$ statements are related to the adjoint action of $G$ on its Lie algebra $\gl$. The proofs of the implications from a collection of statements about the root datum to a statement about the adjoint action are mostly based on known proofs in the case that $\F$ has characteristic $0$. The proofs of the implications from a statement about the adjoint action to a collection of statements about the root datum are different. In this case we assume that one of the statements about the root datum does not hold and then show that the statement about the adjoint action does not hold. For example, we will show that $\neg(7)$ is a consequence of $\neg(1)$ or $\neg(3)$ or $\neg(4)$. The strategy is to make a surjective function from a part of the regular nilpotent elements of the Lie algebra to $\F/\F^{(p)}$ or $\F^\times/(\F^\times)^p$, which is $G$-invariant. For example, in $SL_2(\F)$ with char $\F=2$ we take the function
\[\left(\begin{array}{cc} 0&x\\0&0\end{array}\right)\mapsto x \mod (\F^\times)^2.\]
The proof of $\neg(1)$ or $\neg(3)$ implies $\neg(8)$ is based on the existence of such functions.\\

For $\omega\subset \gl$ define $J(\omega)$ to be the set of distributions with support contained in the closure of ${^G\omega}$. For $L$ a $\co$-lattice in $\gl$ define $J_{L}(\omega)$ to be the image of $J(\omega)$ in the distributions of $\gl/L$ under the canonical map $\phi_L:\gl\rightarrow \gl/L$.

\begin{ver}[Howe]
 For all compact $\omega\subset \gl$ and all $\co$-lattices $L$ of $\gl$:
 \[ \dim J_L(\omega) <\infty.\]
\end{ver}
This conjecture has been proved by Howe in \cite{Ho74} for $G = GL_n(\F)$. Later it has been proved by Harish-Chandra, see \cite{HC99}, for general $G$ in the case that $\text{char}\;\F=0$. In this paper we determine the $\F$-split groups for which Howe's conjecture holds.

Rather surprisingly, these are not only the groups with finitely many nilpotent orbits. Probably there are only finitely many nilpotent orbits with a non-empty intersection with every neighborhood of $0$. To prove Howe's conjecture for certain groups, we will just adapt the proof in \cite{HC99}.

The proof in \cite{HC99} of the local summability of the character of an admissible representation and the local upper bound $|D(g)|^{-\frac12}$ depends on Howe's conjecture. Howe's conjecture is used to proof that the character of a representation is locally a linear sum of Fourier transforms of nilpotent orbital integrals. 

In the following table we list the properties of some $\F$-split groups. In the column nHwC are the set of primes $P$ such that Howe's conjecture does not hold for $G(\F)$ if and only if $\text{char }\F \in P$. In the column INO are the set of primes $P$ such that the group $G(\F)$ has infinitely many nilpotent orbits if and only if $\text{char }\F \in P$.
\[
 \begin{array}{c||c|c|c|c|c}
  G & \text{bad }p & \kappa_v(G) & \rho_v(G)& nHwC & INO\\\hline
  GL_n & - & 1 & 1 & - & -\\
  SL_n & - & n & 1 & p|n & p|n\\
  PGL_n & - & 1 & n & - & p|n\\
  SO_{2n+1} & 2 & 1 & 2 & 2 & 2\\
  SO_{2n} & 2 & 2 & 2 & 2 & 2\\
  Sp_{2n} & 2 & 2 & 1 & 2 & 2\\
  F_4 & 2,3 & 1 & 1 & 2,3 & 2,3\\
  G_2 & 2,3 & 1 & 1 & 2,3 & 2,3\\
  E_8 & 2,3,5 & 1 & 1 & 2,3,5 & 2,3,5
 \end{array}
\]

The obvious direction for generalizing the theory about Howe's conjecture and on the (in)finiteness of nilpotent orbits of this article is to look at reductive groups that are not $\F$-split. The proofs of this article depend heavily on the case by case consideration of the irreducible root systems. It would be nice to find unified proofs.

\section{Notations}
Unless otherwise stated, $\F$ is a local non-Archimedean field with uniformizer $\pi$ and ring of integers $\co$. We define $p := \text{char }\F$. For $n\in \N$ we define $\F^{(n)} := \{x^n : x\in \F\}$ and $\co^{(n)} := \{ x^n : x\in \co\}$.\\

A prime number $p$ is bad for a root system $R$ if
\begin{enumerate*}
 \item $p=2$ and $R$ has a component not of type $A$.
 \item $p=3$ and $R$ has a component of type $E,F$ or $G$.
 \item $p=5$ and $R$ has a component of type $E_8$.
\end{enumerate*}
A prime number $p$ is good for $R$ if it is not bad. See \cite[\S4.1]{SS70} for equivalent definitions of good primes.\\
A prime number $p$ is very good for $R$ if it is good and $R$ does not have a component of type $A_{n}$ with $p$ a divisor of $n+1$.\\
A prime number $p$ is (very) good for $G$ if it is (very) good for the root system of $G$.\\

A $G$-orbit $Ad(G)x$ in $\gl$ is called separable if one of the following equivalent conditions hold
\begin{enumerate*}
\item The differential of the map $g\mapsto Ad(g)x$ is surjective
\item \[\dim \{ g\in G \mid Ad(g)x=x\} = \dim \{ y\in\gl \mid [y,x]=0\} \]
\item The Lie algebra of $\{ g\in G \mid Ad(g)x=x\}$ is equal to $\{ y\in\gl \mid [y,x]=0\}$
\end{enumerate*}

\subsection{$\kappa_v(G)$ \& $\rho_v(G)$}
Let $T$ be a maximal torus of $G$. The embeddings $R(G,T)\hookrightarrow X^*(T)$ and\\ $R^\vee(G,T) \hookrightarrow X_*(T)$ induce group homomorphisms $\Phi : X_*(T)\rightarrow \text{Hom}_\Z(\Z R(G,T),\Z)$ and $\Phi^\vee : X^*(T)\rightarrow \text{Hom}_\Z(\Z R(G,T)^\vee,\Z)$.

\begin{lem}
 Let $\mathscr{T}$ and $\mathscr{S}$ be two complex tori and $\phi : \mathscr{T}\rightarrow \mathscr{S}$. Let $\phi^* : X^*(\mathscr{S})\rightarrow X^*(\mathscr{T})$ be the map $\eps \mapsto \eps\circ \phi$. Then
 \[|\ker \phi / \ker \phi^o| = |(\coker \phi^*)_{tor}|, \]
 where $(\coker\phi^*)_{tor}$ is the torsion part of the cokernel of $\phi^*$.
\end{lem}
\begin{proof}
Choose the bases $\delta_1,\ldots,\delta_m$ for $X^*(\mathscr{S})$ and $\eps_1,\ldots,\eps_n$ for $X^*(\mathscr{T})$ in such a way that
\begin{align*}
 \phi^*(\delta_1) &= d_1\eps_1\\
 \vdots \quad&= \quad\vdots\\
 \phi^*(\delta_k) &= d_k\eps_k\\
 \phi^*(\delta_{k+1},\ldots,\delta_m)&\in \left< d_1\eps_1,\ldots,d_k\eps_k\right>
\end{align*}
then $\prod_{i=1}^k d_i = |(\coker \phi^*)_t|$.\\
Thus 
\begin{align*}
\ker \phi &:= \{ t\in \mathscr{T} \mid \eps_i(t)^{d_i}=1 \text{ for all } 1\leq i\leq k\}\\
\ker \phi^o &:= \{ t\in \mathscr{T} \mid \eps_i(t)=1 \text{ for all } 1\leq i\leq k\}.
\end{align*}
Therefore
\[|\ker \phi/\ker \phi^o| = |\prod_{i=1}^k \Z/d_i\Z|=\prod_{i=1}^k d_i = |(\coker \phi^*)_{tor}|.\qedhere\]
\end{proof}
\begin{gev}
For a complex reductive group, $\mathscr{G}$, $|\coker \Phi|= |\pi_0(Z(\mathscr{G}))|$ and $|\coker \Phi^\vee| = |\pi_1(\mathscr{G}_{\text{der}})|$. 
\end{gev}
\begin{proof}
Look at the adjoint map: $Ad : \mathscr{G}\rightarrow \mathscr{G}^{ad}$. Let $\mathscr{T}$ be a maximal torus of $\mathscr{G}$ and $\mathscr{T}^{ad} = Ad(\mathscr{T})$. Then $\text{Hom}_\Z(\Z\Delta,\Z)=X_*(\mathscr{T}^{ad})$ and $\Phi$ is the map corresponding with $Ad : \mathscr{T}\rightarrow \mathscr{T}^{ad}$:
\[\Phi(\eps):= Ad\circ \eps.\]
We define $\Phi^{tr} : X^*(\mathscr{T}^{ad})\rightarrow X^*(\mathscr{T})$ as follows:
\[ \Phi^{tr}(\eps) := \eps\circ Ad.\]
The cokernel of $\Phi^{tr}$ has a torsion group of order $|\coker \Phi|$. Because $\mathscr{T} \cap \ker Ad = Z(\mathscr{G})$, \[|\coker \Phi|= |Z(\mathscr{G})/Z(\mathscr{G})^o|=|\pi_0(Z(\mathscr{G}))|.\]
Let $\mathscr{G}_{sc}$ be the simply connected cover of $\mathscr{G}_{\text{der}}$. Let $\pi : \mathscr{G}_{sc}\rightarrow \mathscr{G}$ be the following morphism: $\mathscr{G}_{sc}\twoheadrightarrow \mathscr{G}_{\text{der}}\hookrightarrow \mathscr{G}$.\\
Let $\mathscr{T}_{sc}$ be the maximal torus of $\mathscr{G}_{sc}$ such that $\pi(\mathscr{T}_{sc})=\mathscr{T}\cap \mathscr{G}_{\text{der}}$. Then $\text{Hom}_\Z(\Z\Delta,\Z)=X_*(\mathscr{T}_{sc})$ and $\Phi^\vee$ is the map corresponding with $\pi : \mathscr{T}_{sc} \rightarrow \mathscr{T}$. Thus
\[ |\coker \Phi^\vee| = |\ker \pi| = |\pi_1(\mathscr{G}_{\text{der}})|. \qedhere\]
\end{proof}

We call $\rho_v(G) := |\coker \Phi^\vee|$ the virtual order of $\pi_1(G_{\text{der}})$. We call $\kappa_v(G) := |\coker \Phi|$ the virtual number of components of $Z(G)$.

\subsection{Chevalley basis}
The first part of this subsection is based on \cite[\S 1.2]{Ad98}.\\
Let $G$ be a $\F$-split reductive group and $T$ a maximal torus. Let $\gl$ be the Lie algebra of $G$. Let $R := R(G,T)$ be the roots of $G$ and $T$. Let $R^+$ be a set of positive roots of $R$ and $\Delta$ be the set of corresponding simple roots.

We have for $\beta\in R$, the elements $H_\beta$ and $E_\beta$ in $\gl$, such that for all $\alpha,\beta\in R$:
\begin{align*}
 [H_\alpha,H_\beta]&=0\\
 [H_\alpha,E_\beta] &= \left<\alpha^\vee,\beta\right>E_\beta\\
 [E_\beta,E_\alpha] &= \left\{\begin{array}{cl}
                              N_{\beta,\alpha}E_{\beta+\alpha} & \text{if } \beta+\alpha\in R\\
                              H_{\beta} & \text{if } \alpha=-\beta\\
                              0 &\text{otherwise,}
                             \end{array}\right.
 \end{align*}
where each $N_{\beta,\alpha}\in\Z$. For each $\beta\in R$ there exists a unique map $u_\beta : \F\rightarrow G$, such that $d\beta(1)=E_\beta$ and for all $t\in T$ and $x\in\F$, $tu_\beta(x)t^{-1}=u_\beta(\beta(t)x)$. Then $\beta^\vee$, the coroot of $\beta$, is equal to
\[ \beta^\vee(\lambda) = u_\beta(\lambda)u_{-\beta}(-\lambda^{-1})u_\beta(\lambda)u_\beta(-1)u_{-\beta}(1)u_\beta(1).\]
Moreover $d\beta^\vee(1) = H_\beta$.\\
The set $\{H_\alpha : \alpha\in \Delta\}\cup \{ E_\beta : \beta\in R\}$ is called a Chevalley basis. (The term ``basis'' is misplaced here, since if $G$ is not semi-simple it does not span $\gl$ and if $G=PGL_n$ and {char $\F\;|\;n$} it is not linearly independent, see Lemma \ref{ldhpgln}. However if the characteristic is $0$ it is a basis for $\gl'$, the Lie algebra of $G_{\text{der}}=(G,G)$. The $E_\beta$ are always linearly independent.)\\
The adjoint representation $Ad : G\rightarrow End(\gl)$ is determined by the following formulas
\begin{align*}
 Ad(u_\beta(\lambda))E_\alpha &= \left\{ \begin{array}{cl}
                                          E_\beta & \text{if } \beta=\alpha\\
                                          E_{-\beta} + \lambda H_\beta-\lambda^2E_\beta & \text{if } \alpha=-\beta\\
                                          \sum_{i\geq 0} M_{\beta,\alpha,i}\lambda^iE_{i\beta+\alpha} &\text{otherwise}
                                         \end{array}\right.\\
Ad(t)E_\beta &= \beta(t)E_\beta\\
Ad(u_\beta(\lambda))H &= H-d\beta(H)\lambda E_\beta\\
Ad(t)H &= H
\end{align*}
for all $H\in \tl$, the Lie algebra of $T$ and constants $M_{\beta,\alpha,i}\in \F$.\\
The $\F$-points of the image of the algebraic map $Ad$ will be denoted by $Ad(G)$ or $G^{ad}$. From now on we fix a Chevalley basis on $\gl$.

\section{Regular nilpotent orbits}
In the first part of this short introduction to nilpotent orbits, especially regular nilpotent orbits, we will follow \cite[\S 5.1]{Ca85}. Although \cite[\S 5.1]{Ca85} treats regular unipotent elements, we can easily adapt it to regular nilpotent elements.

For each $\alpha\in R$, define $\gl_\alpha := \{ x\in \gl : ad(t)x=\alpha(t)x\}$. We define the height function $ht : R\rightarrow \Z$ as follows:
\[ ht(\sum_{\alpha\in \Delta} c_\alpha \alpha) := \sum_{\alpha\in \Delta} c_\alpha.\]
For $z\in \Z$ we define the following subspaces of $\gl$:
\begin{align*}
\nl_z &:= \bigoplus_{\alpha\in R | ht(\alpha)=z} \gl_\alpha\\
\nl_{\geq z} & := \bigoplus_{\alpha\in R | ht(\alpha)\geq z} \gl_\alpha.
\end{align*}

A nilpotent element of $n\in \gl$ is called a regular nilpotent element if and only if
\[\dim Z_G(n) = \dim T.\]
\begin{pro}\label{Carpros}
Let $\F$ be an algebraically closed field. Let $G$ be a connected reductive group. Then there exist regular nilpotent elements in $\gl$ and any two are conjugated. Let $n\in \gl$ be nilpotent. The following conditions on $n$ are equivalent.
 \begin{enumerate}[label=(\alph*)]
  \item $n$ is regular.
  \item there is a unique Borel subgroup $B$ of $G$ such that $n$ is in the Lie algebra of $B$.
  \item $n$ is conjugated to an element of the form $\sum_{\alpha\in R^+} \lambda_\alpha E_\alpha$ with $\lambda_\alpha \not=0$ for all $\alpha\in \Delta$.
 \end{enumerate}
\end{pro}
\begin{proof}
We use the proof of \cite[Proposition 5.1.2 \& 5.1.3]{Ca85}. That there are only finitely many nilpotent classes is proven in \cite[Theorem 1]{HS85}. The $U$-orbit of $n$ is closed, since every orbit of a unipotent group is closed \cite[Proposition 2.5]{St74}.
\end{proof}
\begin{gev}\label{rnec}
 Let $n,n'$ be regular nilpotent elements of the Lie algebra of $B$. If $g\in G$ is such that $gng^{-1}=n'$, then $g\in B$.\\
If  $n=\sum_{\alpha\in\Delta} c_\alpha E_\alpha$ and $n'=\sum_{\alpha\in\Delta}d_\alpha E_\alpha$, then the following statements are equivalent:
\begin{enumerate*}
\item $n$ and $n'$ are conjugated by an element of $G(\F)$
\item there is a $t\in T$ such that $d_\alpha=\alpha(t)c_\alpha$.
\end{enumerate*}
\end{gev}
\begin{proof}
By Proposition \ref{Carpros} $B=gBg^{-1}$, since $n'$ is in the Lie algebra of $B$ and $gBg^{-1}$. Thus $g\in N_G(B)=B$.\\
Assume that $n=\sum_{\alpha\in\Delta} c_\alpha E_\alpha$ and $n'=\sum_{\alpha\in\Delta}d_\alpha E_\alpha$.\\
If $n$ and $n'$ are conjugated, then there exist  $t\in T$ and $u\in U$ such that $Ad(tu)n=n'$.\\
Since $U$ acts trivial on $\nl/\nl_{\geq 1}$ and $Ad(t)E_\alpha = \alpha(t)E_\alpha$, the second statement follows.\\
If $d_\alpha = \alpha(t)c_\alpha$, then $Ad(t)n=n'$.
\end{proof}
Corollary \ref{rnec} shows that ${^Gn}\cap B={^Bn}$ for all regular $n\in \bl$.

Define $\Phi$ as follows:
\begin{align*}
\Phi: X_*(T)&\rightarrow \hom_\Z(\Z R(G,T),\Z)\\
\Phi:\gamma &\mapsto (\alpha \mapsto \left<\gamma,\alpha\right>)
\end{align*}
The first reason for defining $\Phi$ is the following Proposition. Recall $\kappa_v(G) := |\coker \Phi|$ is the virtual number of components of $Z(G)$.

\begin{pro}\label{ubaan}
 If $G$ is $\F$-split and $p|\kappa_v(G)$, then there are infinitely many regular nilpotent orbits in $\gl$.
\end{pro}
\begin{proof}
Let $\Delta = \{ \alpha_1,\ldots,\alpha_n\}$ be a basis for $R(G,T)$. Define for $1\leq i\leq n$ the function $\epsilon_i \in \hom_\Z(\Z R(G,T) ,\Z)$ by:
\[\epsilon_i(\alpha_j) := \delta_{ij}.\]
So $\epsilon_1,\ldots,\epsilon_n$ is a basis for $\hom_\Z(\Z R(G,T) ,\Z)$. Let $L$ be the image of $\Phi$. Take a compatible basis for $\hom_\Z(\Z R(G,T) ,\Z)$ and $L$: $b_1,\ldots,b_n$ and $d_1b_1,\ldots,d_nb_n$ with $d_i|d_{i+1}$. Since the cokernel is finite it has $\prod_{i=1}^n d_i$ elements. Define $M \in GL_n(\Z)$ by
\[ \left(\begin{array}{c} b_1\\ \vdots \\ b_n\end{array}\right) = M \left(\begin{array}{c} \epsilon_1\\ \vdots \\ \epsilon_n\end{array}\right).\]
Look at the following subset of $\nl_1$:
\[\nl'_1 := \{ \sum_{\alpha\in \Delta} c_\alpha E_{\alpha} : c_\alpha\in\F^\times \}.\]
Define $\pi$ to be the following parametrization of $\nl'_1$:
\[ \pi :(\F^\times)^n \rightarrow \nl'_1,\quad \pi(c_1,\cdots,c_n) := \sum_{i=1}^n c_iE_{\alpha_i}.\]
Since $\sum_{\alpha\in\Delta}c_\alpha E_\alpha$ is regular, it is in the same conjugacy class of $G(\F)$ as 
$\sum_{\alpha\in\Delta}d_\alpha E_\alpha$ if and only if there is a $t\in T$ such that $d_\alpha = \alpha(t)c_\alpha$ for all $\alpha\in \Delta$ by Corollary \ref{rnec}.\\
Let $A\in GL_n(\Z)$ and define $\phi_A:(\F^\times)^n\rightarrow (\F^\times)^n$ by:
\[\phi_A(x_1,\ldots,x_n) := (\prod_{i=1}^n x_i^{a_{1i}} , \ldots, \prod_{i=1}^n x_i^{a_{ni}}).\]

Now $\pi\circ \phi_A$ is also a parametrization of $\nl'_1$ and
$$\phi_A\pi^{-1}(t\pi(\phi_{A^{-1}}(x_1,\ldots,x_n))t^{-1}) = ((\prod_{i=1}^n\alpha_i(t)^{a_{1i}})x_1,\ldots,(\prod_{i=1}^n\alpha_i(t)^{a_{ni}}) x_n ).$$
Define this action of $T$ on $(\F^\times)^n$ to be the action with respect to $A$.\\
Take $A := (M^{-1})^t$.\\
We claim that for every $\gamma\in X_*(T)$ with $\Phi(\gamma) = \sum_{i=1}^n z_id_ib_i$ one has the following action on $(\F^\times)^n$ with respect to $A$ of $\gamma(s)$:
$$(x_1,\ldots, x_n)\mapsto (s^{z_1d_1} x_1,\ldots,s^{z_nd_n}x_n).$$
To prove this claim, consider the factor in front of $x_j$:
\begin{align*}
 \prod_{i=1}^n \alpha_i(\gamma(s))^{a_{ji}}&=s^{\sum_{i=1}^n a_{ji}\left<\gamma,\alpha_i\right>}
\intertext{Evaluate the power of $s$:}
 \sum_{i=1}^n a_{ji}\left<\gamma,\alpha_i\right>&=\sum_{i=1}^n a_{ji}\sum_{k=1}^n z_kd_kb_k(\alpha_i) = \sum_{i=1}^n a_{ji} \sum_{k=1}^n z_kd_k m_{ki} = \sum_{k=1}^n z_kd_k \sum_{i=1}^n a_{ji}m_{ki}.\\
\intertext{Since $A = (M^{-1})^t$ one has that $\sum_{i=1}^n a_{ji}m_{ki}=\delta_{jk}$, therefore}
\sum_{i=1}^n a_{ji}\left<\gamma,\alpha_i\right> &= z_jd_j\text{ , hence}\\
\prod_{i=1}^n \alpha_i(\gamma(s))^{a_{ji}}&=s^{z_jd_j}.
\end{align*}

Since $p | \# \coker \Phi$, then $p|d_n$. Identify $\nl'_1$ with $(\F^\times)^n$ via the parametrization $\pi \circ \phi_{M^t}$. Look at the $n$-th coordinate: $x_n\mapsto s^{z_\gamma d_n}x_n$ for every pair $\gamma \in X_*(T)$, $s\in \F^\times$. The images of the cocharacters generate the torus, so the orbit of the $n$-th coordinate under $T$ is contained $\{s^{d_n}x_n : s\in \F^\times \}$. Hence if $(x_1,\ldots,x_n)$ is in the same orbit as $(y_1,\ldots,y_n)$, then there is a $s\in \F^\times$ such that $s^{d_n}x_n=y_n$. Because $p|d_n$ the group $\F^\times/(\F^\times)^{d_n}$ is infinite. We conclude that if $p$ divides the order of the cokernel, then there are infinitely many regular nilpotent orbits.
\end{proof}

\section{The virtual number of components of $Z(G)$ and ...}
As we saw in the previous section, when $p$ divides the virtual number of components of $Z(G)$ there are infinitely many regular nilpotent orbits. In this section we show that even more properties that hold when the characteristic is zero, do not hold anymore when $p|\kappa_v(G)$. After giving these counter examples for theorems that hold in characteristic zero, at the end of this section we show that $p|\kappa_v(G)$ for a restrictive class of reductive groups. By the way, the condition $p|\kappa_v(G)$ is based on the group $SL_n(\F)$ with $p|n$. As turns out in the end for $p\geq 5$ the main example is $SL_n(\F)$. In this section all the properties are geometric in nature, so we do not have to worry about rationality.

\subsection{Separability and $\kappa_v(G)$}
\begin{lem}
Let $X := \sum_{\alpha\in\Delta}E_\alpha$. The map $[X,\cdot]:\tl\rightarrow \nl_1$ is not surjective if and only if $p|\kappa_v(G)$.
\end{lem}
\begin{proof}
See $\tl$ as $X_*(T)\otimes \F$. Let $Y\in X_*(T)$, then $[dY(1),E_\alpha]=\left<\alpha,Y\right>E_\alpha$. Let $X_1,\ldots,X_n$ be a basis for $X_*(T)$. The matrix $M$ corresponding to $\Phi$ with respect to the basis $X_1,\ldots,X_n$ and the dual basis of $\Delta$ in $\text{Hom}_\Z(\Z \Delta,\Z)$ is the same as the matrix corresponding to $[X,\cdot]$ with respect to the basis $X_1,\ldots,X_n$ and $(E_\alpha : \alpha\in \Delta)$. Let $d_1,\ldots,d_{|\Delta|}$ be the integers on the diagonal of the Smith normal form of $M$. Then $\kappa_v(G) =\# \coker \Phi = \prod_{i=1}^{|\Delta|} d_i$. Also there are $E_1,\ldots, E_{|\Delta|}$ such that $\nl_1 = \left<E_1,\ldots,E_{|\Delta|}\right>$ and $[X,\tl] = \left<d_1E_1,\ldots,d_nE_{|\Delta|}\right>$.
\end{proof}
Recall an $G$-orbit $Ad(G)x$ in $\gl$ is separable if and only if
\[\dim \{ g\in G \mid Ad(g)x=x\} = \dim \{ y\in\gl \mid [y,x]=0\}.\]
\begin{ste}\label{caisrn}
 If $p|\kappa_v(G)$, then the regular orbit is not separable.
\end{ste}
\begin{proof}
Let $T$ be a torus and $X\in \nl_1$ a regular element. Then $\dim Z_T(X) = \dim T-|\Delta|$ by Proposition \ref{Carpros}. Thus if the orbit of $X$ is separable, then $[X,\cdot] : \tl \rightarrow \nl_1$ has a kernel of dimension $\dim T-|\Delta|$. Thus $[X,\cdot]$ must be surjective. Since $p|\kappa_v(G)$ the map $[X,\cdot]$ is not surjective. Hence the orbit of $X$ is not separable.
\end{proof}

\subsection{ $Ad$ and $\kappa_v(G)$}
Let $G$ be a reductive $\F$-group. We will go back and forth between $G$ and $Ad(G)$. Therefore we have a look at the adjoint map $Ad : G\rightarrow Ad(G)$. The adjoint map is defined over $\F$. We will show that $d(Ad)$ maps non-zero-nilpotent elements to non-zero nilpotent elements. $Ad$ is separable (ie, $d(Ad)$ is surjective) if and only if $p{\not|}\kappa_v(G)$. To distinguish the objects associated with $Ad(G)$ from the ones associated with $G$, the ones associated with $Ad(G)$ get a superscript $ad$: $G^{ad}$,$\gl^{ad}$, $\nl^{ad}$, ect.
\begin{lem}\label{darnib}
 $d(Ad) : \nl\rightarrow \nl^{ad}$ is an isomorphism.
\end{lem}
\begin{proof}
Take a Chevalley basis on $\gl$. Let $\alpha\in R(G,T)$. Let $G^{ad}_\alpha$ be the image of $Ad\circ u_\alpha : \F\rightarrow G^{ad}$. The action of $u_\alpha(x)$ on certain elements of $\gl$ is as follows: 
\begin{align*}
Ad(u_\alpha(x))E_{-\alpha} &= E_{-\alpha}+xd\alpha^\vee(1)-x^2E_\alpha,\\
Ad(u_\alpha(x)) H&=H-d\alpha(H)xE_\alpha.
\end{align*}
Since $\left<\alpha,\alpha^\vee\right>=2$, either $d\alpha^\vee(1) \not=0$ or there exists a $H\in\tl$ such that $d\alpha(H)\not=0$. Therefore $Ad\circ u_\alpha$ is an isomorphism between $\F$ and its image in $G$. Because $tu_\alpha(x)t^{-1}=u_\alpha(\alpha(t)x)$ for $t\in T$ and $x\in \F$, also
\[tAd(u_\alpha(x))t^{-1} = Ad(u_\alpha(\alpha(t)x)), \]
for all $t\in T^{ad}$ and $x\in \F$. Thus $d(Ad) : \gl_\alpha\rightarrow \gl_\alpha^{ad}$ is an isomorphism. Therefore $d(Ad) : \nl\rightarrow \nl^{ad}$ is injective. Since $\dim \nl=\dim \nl^{ad}$, the Lemma follows.
\end{proof}

\begin{pro}\label{csad}
 The map $Ad : G\rightarrow Ad(G)$ is separable if and only if the characteristic of $\F$ does not divide the virtual number of components of $Z(G)$.
\end{pro}
\begin{proof}
Let $\Delta$ be a system of positive roots for $R(G,T)$.\\
Define $n := |\Delta|$ and let $\alpha_1,\ldots,\alpha_n$ be the roots in $\Delta$. Take $\gamma_1,\ldots,\gamma_n\in X_*(T)$ such that the image of $\Phi$ is generated by $\gamma_1,\ldots,\gamma_n$. The number of elements in the cokernel of $\Phi$ is equal to the determinant of the matrix $M_{ij} := \left<\gamma_j,\alpha_i\right>$. Since $d(Ad)$ is surjective on $\nl_+^{ad}\oplus \nl_-^{ad}$, we only have to look whether $Ad : T\rightarrow T^{ad}$ is separable. Identify $T^{ad}$ with a torus of dimension $n$ in such a way that the map $Ad$ is as follows:
\[t\mapsto \left(\begin{array}{ccc} \alpha_1(t) & & 0\\ &\ddots&\\ 0& & \alpha_n(t)\end{array}\right) \]
The Lie algebra of a torus $S$ is canonically isomorphic to $X_*(S)\otimes_\Z \F$ \cite[4.4.11(4)]{Sp98}. With this isomorphism the map $d(Ad)$ is the linear map such that for $\gamma\in X_*(T)$, $d(Ad)(\gamma)= Ad\circ \gamma$. Now the images of $Ad \circ \gamma_1,\ldots,Ad\circ\gamma_n$ generate the image of $d(Ad)$. Thus the image of $\tl$ is generated by the vectors $\sum_{i=i}^n \left<\gamma_j,\alpha_i\right>\chi_i$ for $j=1,\ldots,n$. This is surjective if and only if the corresponding matrix has non-zero determinant. The corresponding matrix is equal of $M$. Thus $p{\not|}\kappa_v(G)=\#\coker \Phi$ if and only if $M$ is invertible if and only if $Ad$ is separable.
\end{proof}

\subsection{Very good primes and $\kappa_v(G)$}
\begin{lem}
 If $Y\subset X_*(T)$ such that $\Phi(Y)$ has finite index in $\hom(\Z R(G,T),\Z)$, then $\#\coker \Phi$ divides $\#\coker \Phi|_Y$.
\end{lem}
\begin{proof} The lemma follows from general abstract non-sense:
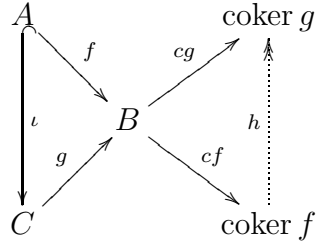
\[
\xymatrix{A \ar[dr]^f \ar@{^{(}->}[dd]^\iota & & \coker g \\ & B\ar[ur]^{cg} \ar[dr]^{cf} & \\ C \ar[ur]^g & & \coker f \ar@{.>>}[uu]^h}
\]
Since $cg\circ g \circ \iota =0$ there is an unique morphism $h : \coker f \rightarrow \coker g$ such that $cg = cf \circ h$. Since $cg$ is surjective, also $h$ is surjective. Thus $\#\coker g |\# \coker f$. \end{proof}

\begin{pro} 
 If $p|\kappa_v(G)$, then $p$ divides the determinant of the Cartan matrix of $R(G,T)$.
\end{pro}
\begin{proof}
Let $Y$ be the subgroup of $X_*(T)$ generated by the coroots of $R(G,T)$. The order of the cokernel $Y\rightarrow \hom_Z(\Z R(G,T),\Z)$ is equal to the determinant of the Cartan matrix.
\end{proof}
\begin{gev}\label{2implies134}
 If $p|\kappa_v(G)$, then $p$ is not a very good prime for $G$. If moreover $G$ does not contain a normal subgroup of type $A_l$, then $p$ is a bad prime for $G$ and $p\in\{2,3\}$.
\end{gev}
\begin{proof}
By \cite[11.4, Exercise 2]{Hu78} the determinants of the Cartan matrices for the irreducible root systems are:
$$A_l: l+1;B_l:2;C_l:2;D_l:4; E_6:3; E_7:2; E_8,F_4\text{ and } G_2: 1.$$
Compare this with the notion of a prime that is not a very good prime. Then $p$ divides:
\[A_l: l+1;B_l:2;C_l:2;D_l:2; E_6:2,3; E_7:2,3; E_8: 2,3,5; F_4:2,3; G_2: 2,3.\qedhere\]
\end{proof}

\section{Howe's conjecture in bad characteristic}
In this section we show that Howe's conjecture does not hold for $\F$-split groups in bad characteristic. The calculations in the actual group are postponed to the end of this section and the Appendix. Under the assumption that there exists a bad pair, we will construct sets of linearly independent distributions in $J_L(\omega)$ of arbitrary finite size. The support of these distributions is contained in the set of nilpotent elements. Two consequences of our method are the existence of infinitely many regular nilpotent orbits and the inseparability of the regular nilpotent orbit.

\subsection{Reduction to bad pairs}\label{redbp}
Let $G$ be a $\F$-split reductive group. Let $T$ be a maximal $\F$-split torus. Let $R^+$ be a system of positive roots. Let $U^+$ be the unipotent subgroup corresponding to $R^+$ and $\nl$ its Lie algebra. Let $B=TU^+$ be the corresponding Borel subgroup. The set of regular nilpotent elements of $\nl$ is denoted by $\nl'$.

Let $H_1,\ldots,H_r$ and $E_\gamma$ for $\gamma\in R$ be a Chevalley basis for $\gl$. Let $u_\gamma : \F\rightarrow U_\gamma$ be the corresponding parametrization of $U_\gamma$. ($du_\gamma(1)=E_\gamma$) Now $\nl$ has as basis $E_\alpha : \alpha\in R^+$.

\begin{mdef}
Let $\eta : \F \rightarrow \nl'$ and $\chi : \nl'\rightarrow \F$ be polynomial functions. The pair $(\eta,\chi)$ is called a bad pair if it satisfies the following four conditions:
\begin{enumerate*}
 \item $\chi \eta(\alpha)=\alpha$ for all $\alpha\in\F$.
 \item If $n,n'\in \nl'$ are conjugated, then $\chi(n)\equiv \chi(n') \mod \F^{(p)}$. There exists a $z\in p\Z$ such that $c^z\chi(n) = \chi(cn)$. \label{Binv2}
 \item For $\gamma\in R^+$ and $\alpha\in \co^\times$, $\eta(\alpha)_{\gamma} \in \co$ and if moreover $\gamma\in \Delta$, then $\eta(\alpha)_{\gamma} \in \co^\times$.
 \item $\chi \in \co[X_\gamma,X^{-1}_\beta : \gamma\in R^+, \beta\in \Delta]$.
\end{enumerate*}
\end{mdef}
For the remainder of this subsection we assume that $(\eta,\chi)$ is a bad pair. Because $\F/\F^{(p)}$ is infinite, the first and second conditions of a bad pair already imply that there are infinitely many regular nilpotent orbits in $\gl$. We will use $\chi$ to define $G$-invariant distributions and $\eta$ to show that they are linearly independent.\\

For $n\in \N$ we define 
\begin{align*}
 U_{\gamma,n} &:= u_\gamma(v^{-1}[n,\infty)),\\
 T_{i,n} & := \{ t\in T \mid \forall[\alpha\in X^*(T)]\; v(\alpha(t)-1)\geq n\}.
\end{align*}
For $n\in \N$ define the group $K_n$ to be the group generated by the groups $U_{\gamma,n}$ and $T_{i,n}$. Define $K := K_0$.\\
We may identify $\tl$ with $X_*(T)\otimes_\Z\F$ by 
\[ X_*(T)\ni \gamma \mapsto d\gamma(1) \in \tl.\]
Let $\delta_1,\ldots,\delta_s$ a basis for $X_*(T)$ and $H'_1,\ldots,H'_s$ the corresponding basis in $\tl$.\\
Let $L$ be the $\co$-lattice spanned by $H'_1,\ldots,H'_s$ and all $E_\gamma$.\\
For $m\in \gl$ we define $m_i\in\F$ and $m_\gamma\in\F$ such that
\[ m = \sum_{i=1}^s m_iH'_i + \sum_{\gamma\in R} m_\gamma E_\gamma.\]
Now $L$ is $K$-invariant. Thus $K$ acts on $L/\pi^nL$. The group $K_n$ acts trivial on $L/\pi^n L$, by the choice of $K_n$ and $L$. Now $K/K_1\cong G(\F_q)$ and $\gl(\F_q)=L/\pi L$.

\begin{lem}\label{kiesN3}
 There exists a $N>0$ such that for all $n\in\N_{>0}$, $k\in K$ and $\alpha\in \co^\times$:
\[k\eta(\alpha)k^{-1} \in \nl + \pi^{Nn}L \Rightarrow k\in (B\cap K)K_n.\]
\end{lem}
\begin{proof}
 The map $\pi_0: G(\co)\rightarrow G(\F_q)$ gives a corresponding map on the Lie algebra: $\pi_0: \gl(\co) \rightarrow \gl(\F_q)$ with kernel $\pi L$. Since $\pi_o(\eta(\alpha))$ is also a regular nilpotent element and $\pi_0(k\eta(\alpha)k^{-1})\in \nl(\F_q)$, we have $\pi_0(k)\in B(\F_q)$. Thus $k\in (B\cap K)K_1$. Take for the moment a general $N\in \N_{>0}$. Because $\mathfrak{n}$ and $\pi^{Nn}L$ are $(B\cap K)$-invariant and $K_1 = (B\cap K_1)(U^-\cap K_1)$, we may assume $k\in U^-\cap K$. Take $x_{\gamma}\in \co$ for $\gamma\in R^-$ such that $u = \prod_{\gamma\in R^-}u_\gamma(x_\gamma)$. Let $p_{i},p_\beta \in \F[X_\gamma:\gamma\in R^{-},Y,Y^{-1}]$ be such that $(u\eta(\alpha)u^{-1})_{i} = p_{i}(x_\gamma,\alpha,\alpha^{-1})$ and $(u\eta(\alpha)u^{-1})_{\beta} = p_{\beta}(x_\gamma,\alpha,\alpha^{-1})$.\\
 Let $I$ be the ideal generated by $p_{\beta}$ for $\beta\in R^-$. Then $u\eta(\alpha)u^{-1}\in \nl$ if and only if $p_{\beta}(x_\gamma,\alpha,\alpha^{-1})=0$ for all $\beta\in R^-$. Because of Corollary \ref{rnec} for $x_\gamma,\alpha\in\overline{\F}$: 
 \[u\eta(\alpha)u^{-1}\in \nl \Leftrightarrow u=1 \Leftrightarrow \forall\gamma\in R^{-}[x_\gamma=0],\]
 where $u=\prod_{\gamma\in R^-} u_\gamma(x_\gamma)$. By the Nullstellensatz we have $X_\gamma \in \sqrt{I}$ for all $\gamma\in R^-$. Thus there exists a $m\in \N$ such that $X_\gamma^m \in I$ for all $\gamma\in R^-$. Therefore there are polynomials $f_{\gamma,\beta}\in \F[X_\gamma,Y,Y^{-1}]$ such that
 \[X_\gamma^m = \sum_{\beta\in R^+} f_{\gamma,\beta}p_{\beta}.\]
 Let $M$ be the smallest $n\in \N_{\geq 0}$ such that $f_{\gamma,\beta}(x_\gamma,\alpha,\alpha^{-1}) \subset \pi^{-n}\co$ for all $\beta,\gamma\in R^-$,  $x_\gamma\in \co$ and $\alpha\in\co^\times$.\\
 Take $N := m+M$. Assume that $u\eta(\alpha)u^{-1} \in \nl+\pi^{Nn}L$, then $v(p_{\beta}(x_\gamma,\alpha,\alpha^{-1}))\geq Nn$. Because 
 \[x^m_\gamma = \sum_{\beta\in R^-}f_{\gamma,\beta}(x_\gamma,\alpha,\alpha^{-1})p_{\beta}(x_\gamma,\alpha),\] 
 we have 
 \[v(x^m_\gamma)\geq Nn-M=mn+(n-1)M\geq mn.\]
 Thus $v(x_\gamma)\geq n$.
\end{proof}

Let $\delta_B$ be the modular function of $B$, thus
\[\delta_B(b)\int_B f(xb)dx=\int_B f(x)dx.\]

\begin{pro}[Rao]\label{disgi2}
 Assume that $V\subset \nl$ is open and $B$-invariant. Then for all $f\in C_c^\infty(\gl)$,
\[\int_V f(bXb^{-1})dX =\delta_B(b)\int_Vf(X)dX.\]
 Moreover the distribution
 \[ D_V(f) := \int_V \int_K f(kXk^{-1})dkdX\]
 is $G$-invariant.
\end{pro}
\begin{proof}
Since $\delta|_B(b) = |\det Ad\; b|_\nl|^{-1}$,
\[ \int_\nl f(bXb^{-1})dX = \delta_B(b)\int_{\nl} f(X)dX.\]
Because $V$ is open and $B$-invariant, we can apply this formula to $\int_V f(X)dX$. This proves the first statement of the Proposition.\\
The second statement follows from the first by \cite[Proposition 4]{Ho74}. The method described here is essentially in \cite{Ra72}.
\end{proof}
\begin{gev}\label{deojw}
 Let $\omega\subset \gl$ be open and compact. If $V\subset \nl$ is open and $B$-invariant, then $D_V\in J(\omega)$.
\end{gev}
\begin{proof}
 By Proposition \ref{disgi2}, $D_{V}$ is a $G$-invariant distribution. The support of $D_{V}$ is contained in $\overline{\nl^K}$. Since $\omega$ is open, there is a $m\in\N$ such that $\pi^m\nl(\co)\subset \omega$. Since $\nl(\co)^T = \nl$, then $\text{supp } D_{V}\subset \overline{\nl^K}\subset \overline{(\omega^T)^K}\subset \overline{\omega^G}$.
\end{proof}

For $\alpha\in \F^\times$ and $s\in \N$, define $V_{\alpha,s}\subset \nl$ as follows:
\[V_{\alpha,s} := \{ n\in\nl \mid \chi(n) \equiv \alpha \mod (\pi^s\co+\F^{(p)})\}.\]
Let $\Delta = \{ \alpha_1,\ldots,\alpha_m\}$. Define for $a_1,\ldots,a_m\in \F$ the following nilpotent element:
\[ n(a_1,\ldots,a_m) := \sum_{i=1}^m a_i E_{\alpha_i}.\]
Take $z\in p\Z$ such that $\chi(cn)=c^z\chi(n)$ for all $n\in\nl'$ and $c\in \F$. 
\begin{lem}\label{a=bmodnp}
Let $n\in\N_{>0}$, $\alpha\in\co^\times$ and $\beta\in \co^\times$.\\
If
 \[ \int_{V_{\pi^{-znN}\beta,n}}\int_{k\in K} 1_{\pi^{-Nn}\eta(\alpha)+L}(kXk^{-1})dkdX>0,\]
then $\alpha \equiv \beta \mod \pi^n\co+\co^{(p)}$.
\end{lem}
\begin{proof}
Let $X\in V_{\pi^{-znN}\beta,n}$, $k'\in K$ and $l'\in L$ such that
\[ k' X k'^{-1}+l' = \pi^{-Nn}\eta(\alpha).\]
Since $L$ is $K$-invariant, there exist a $k\in K$ and $l\in L$ such that
\[k\pi^{-nN}\eta(\alpha)k^{-1}+l = X\in V_{\pi^{n-znN}\beta,n}\subset \nl.\]
Thus $k\in (K\cap B)K_n$ by Lemma \ref{kiesN3}, because $k\eta(\alpha)k^{-1}\in \nl+\pi^{nN}L$. Take $b_k\in K\cap B$ and $k_n\in K_n$ such that $k=k_nb_k$. Take $a_1,\ldots,a_m\in \co^\times$ and $n_2\in \nl_2(\co)$ such that $b_k\eta(\alpha)b_k^{-1}=n(a_1,\ldots,a_m)+n_2$. By assumption \ref{Binv2} of the bad pair, there exists a $\gamma\in\F$ such that $\chi(n(a_1,\ldots,a_m)+n_2) = \alpha+\gamma^p$.  Since $k_n\in K_n$ and $n(a_1,\ldots,a_m)+n_2\in L$, there exists a $l'\in L$ such that\\ $k_n(n(a_1,\ldots,a_m)+n_2)k_n^{-1}= n(a_1,\ldots,a_m)+n_2+\pi^nl'$. Thus
\begin{align*}
 \chi(k\eta(\alpha) k^{-1}+\pi^{nN}l) &=\chi(k_nb_k\eta(\alpha)b_k^{-1}k_n^{-1}+\pi^{nN}l)\\
 &= \chi(n(a_1,\ldots,a_m)+n_2+\pi^n l'+\pi^{nN}l).
\intertext{Since the $a_i$ are in $\co^\times$ and $\chi \in\co[X_\gamma,X^{-1}_\beta: \gamma\in R^+,\beta\in \Delta] $,}
 \chi(n(a_1,\ldots,a_m)+n_2+\pi^n l'+\pi^{nN}l)&\equiv\chi(n(a_1,\ldots,a_m)+n_2)\\&=\alpha+\gamma^p\mod \pi^n\co.
\intertext{Because $\chi(\pi^{-nN}x)=\pi^{-znN}\chi(x)$ for all $x\in \gl$,}
 \chi(k\pi^{-nN}\eta(\alpha)k^{-1}+l) &\equiv (\alpha+\gamma^p)\pi^{-znN} \mod \pi^{n-znN}\co. 
\intertext{
Since $k\pi^{-nN}\eta(\alpha)k^{-1}+l\in V_{\pi^{-znN}\beta,n}$,} \chi(k\pi^{-nN}\eta(\alpha)k^{-1}+l) &\equiv \pi^{-znN}\beta \mod (\pi^{n-znN}\co+\F^{(p)}).\\ 
\intertext{Thus}
 \pi^{-znN}\beta \equiv \chi(k\pi^{-nN}\eta(\alpha)k^{-1}+l)&\equiv \pi^{-znN}\alpha \mod (\pi^{n-znN}\co+\F^{(p)}).
\end{align*}
Then $\alpha \equiv \beta \mod (\pi^n\co+\F^{(p)})$. Because $\F^{(p)}\cap \co = \co^{(p)}$ and $\alpha,\beta\in \co$, the Lemma follows.
\end{proof}

\begin{ste}\label{bpnhc}
 Let $G$ be a $\F$-split reductive group. If there exists a bad pair $(\eta,\chi)$ for $G$, then $\dim J_L(\omega)= \infty$.
\end{ste}
\begin{proof}
 Take $n\in \N_{>0}$. Let $\alpha_1,\ldots,\alpha_k$ be representatives of the cosets of $\pi^n\co+\co^{(p)}$ in $\co$. Define for $1\leq i\leq k$ the following distribution and function:
 \begin{align*}
 D_i(f) &:= D_{V_{\pi^{-znN}\alpha_i,n}}(f) = \int_{V_{\pi^{-znN}\alpha_i,n}}\int_{k\in K} f(kXk^{-1})dkdX,\\
 f_i & := 1_{\pi^{-nN}\eta(\alpha_i)+L}.
 \end{align*}
 The distribution $D_i$ are in $J(\omega)$ by Corollary \ref{deojw}. Let $c_i := D_i(f_i)$, then $c_i>0$. The to $C_c^\infty(\gl/L)$ restricted distributions $D_1,\ldots,D_k$ are linearly independent, since $D_i(f_j) = c_i \delta_{ij}$ by Lemma \ref{a=bmodnp}. Thus $\dim J_L(\omega)\geq k$. As $n$ goes to infinity so does $k$.
\end{proof}

\subsection{The bad pair construction}
In this section we assume that char $\F$ is bad for $G$. The construction of a bad pair is done in three steps. First we construct a bad pair in the case $G$ is simple of adjoint type. Then we show that if there is a bad pair for $Ad(G)$, then we can construct a bad pair for $G$. In the third step we combine the results of the first and second step to construct a bad pair.\\

Define $X := \sum_{\beta\in \Delta} E_\beta$. Let $\alpha_1,\ldots,\alpha_k$ be the roots of height $p+1$. Define 
\[n(a_1,\ldots,a_k) := X+\sum_{i=1}^k a_i E_{\alpha_i},\]
for $a_i\in\F$.
\begin{lem}\label{cbcpl}
 If $p$ is bad for the simple group $G$, then
 \begin{enumerate}
  \item $[X,\nl_i] = \nl_{i+1}$ if $i<p$. 
  \item $\dim \nl_{p+1}/[X,\nl_p] = 1$
  \item $\dim \nl_{1}=\dim \nl_{p+1} +1$
 \end{enumerate}
\end{lem}
\begin{proof}
 Let $U$ be the unipotent subgroup of the Borel subgroup $B$. Let $Ad : G\rightarrow Ad(G)$ be the surjective homomorphism between $G$ and its adjoint representation on the Lie algebra. Then $B^{ad} := Ad(B)$ is a Borel subgroup of $Ad(G)$ and $U^{ad} := Ad(U)= R(B^{ad})$. The map $Ad: U\rightarrow U^{ad}$ is bijective and its tangent map $d(Ad)$ is also bijective by Lemma \ref{darnib}. Therefore $\nl\cong d(Ad)(\nl)=\nl^{ad}$, where $\nl^{ad}$ is the Lie algebra of $U^{ad}$. Since $X$ and $\nl_i$ are all in $\nl$ and $Ad(G)$ is of adjoint type, the Lemma only depends on the root system (and not on the root datum). So we just have to go through the root types. In the Appendix the Lemma is checked for the root data of adjoint type.
\end{proof}
\begin{gev}\label{rnons}
 If $p$ is bad for $G$, then the regular nilpotent orbit is not separable.
\end{gev}
\begin{proof}
Since $X$ is a regular nilpotent element
\[ \dim Z_G(X) = \dim Z_B(X) = \dim B-\dim U= \dim T.\]
Because $\dim Z_T(X) = \dim T-|\Delta|$, we have $\dim Z_U(X) = |\Delta|$.\\
Thus if the orbit of $X$ is separable, then $[X,\cdot]: \nl \mapsto \nl_{\geq2}$ has a kernel of dimension $|\Delta|$. Therefore $[X,\cdot]$ must be surjective. As Lemma \ref{cbcpl} shows, this is not the case when $p$ is bad for a simple group $G$. By passing to the adjoint group, the Corollary follows. 
\end{proof}

\begin{pro}\label{linfun}
 There exists a surjective linear function $f : \F^k\rightarrow \F$ such that
 if $n(a_1,\ldots,a_k)$ is conjugated to $n(b_1,\ldots,b_k)$, then $f(a_1,\ldots,a_k)\equiv f(b_1,\ldots,b_k) \mod \F^{(p)}$.
\end{pro}
\begin{proof}
 Let $f : \nl_{p+1} \rightarrow \F$ be a linear function corresponding with the isomorphism $\nl_{p+1}/[X,\nl_p]\cong \F$. For $u\in U$ write $u = \prod_{\gamma\in R^+}u_\gamma(x_\gamma)$.\\
 By Lemma \ref{cbcpl}, we have $\dim \nl_{i} = \dim \nl_1-1$ for $2\leq i\leq p+1$ and $n\mapsto [X,n]$ is a bijection from $\nl_i$ to $\nl_{i+1}$ for $2\leq i \leq p-1$.\\
 
 We will prove with induction on the height of the roots that there exist $c_\gamma,d_\gamma\in\F$ and $x\in \F$, such that for $i\leq p-1$, if $uXu^{-1} \equiv X\mod \nl_{\geq i+2}$, then 
 \begin{align} x_\gamma &= c_\gamma x^{ht(\gamma)},\label{stingev1}\end{align} for $\gamma\in R^+$ with $ht(\gamma)\leq i-1$ and 
 \begin{align}
 uXu^{-1}&\equiv X-[X,\sum_{\gamma\in R^+_{i+1}} x_\gamma E_\gamma] + \sum_{\gamma\in R^+_{i+2}} d_\gamma x^{i+1} E_\gamma \mod \nl_{\geq i+3}.\label{stingev2}
 \end{align}
 
 Before we give the induction argument, first we restate (\ref{stingev2}).\\
 The nilpotent element $uXu^{-1} \mod \nl_{\geq i+3}$ only depends on the value of $x_\gamma$ for the $\gamma$ with height at most $i+1$. In expression (\ref{stingev2}) the dependence of the roots of height $i+1$ is taken care of with the term $-[X,\sum_{\gamma\in R^+_{i+1}} x_\gamma E_\gamma]$. So for the proof of (\ref{stingev2}) we need to show that
 \[RM_i := uXu^{-1}-X+[X,\sum_{\gamma\in R^+_{i+1}} x_\gamma E_\gamma] \]
 is equal to $\sum_{\gamma\in R^+_{i+2}} d_\gamma x^{i+1} E_\gamma$, when $uXu^{-1} \equiv X \mod \nl_{\geq i+2}$.\\
 
 The function $[X,\cdot]:\nl_1\rightarrow \nl_2$ gives that $x := x_\gamma=x_\delta$ for all $\gamma,\delta\in\Delta$. By the Steinberg conjugacy formula \cite[Proposition 8.2.3]{Sp98} we have $d_\gamma\in\F$ such that $RM_{1} = \sum_{\gamma\in R^+_{3}} d_\gamma x^2 E_\gamma$.\\
 Assume that (\ref{stingev1}) and (\ref{stingev2}) hold for $i-1$.\\
 The function $[X,\cdot]:\nl_i\rightarrow \nl_{i+1}$ is bijective. So for all $n_{i+1}\in\nl_{i+1}$ there is exactly one $u\in U_i$ such that $u(X+n_{i+1})u^{-1} = X \mod \nl_{\geq i+2}$, namely the one corresponding with the inverse of  $[X,\cdot]$. Let $In : \nl_{i+1}\rightarrow \nl_i$ be the inverse of $[X,\cdot]$. Then the $x_\gamma\in\F$ for $\gamma\in R^+_i$ are such that
 \[In(RM_{i-1}) = \sum_{\gamma\in R_i^+}x_\gamma E_\gamma.\]
 By the induction hypotheses $RM_{i-1} = \sum_{\gamma\in R^+_{i+1}} d_\gamma x^{i} E_\gamma$ for some constants $d_\gamma\in \F$. Thus for every $\gamma \in R_i^+$ there exists a $c_\gamma\in \F$ such that $x_\gamma = c_\gamma x^{i}$. By the Steinberg conjugacy formula we have $d_\gamma\in\F$ such that $RM_{i} = \sum_{\gamma\in R^+_{i+2}} d_\gamma x^{i+1} E_\gamma$.\\
 
 Assume that $un(a_1,\ldots,a_k)u^{-1} \equiv n(b_1,\ldots,b_k) \mod \nl_{\geq p+2}$.\\
 Then certainly 
 \[uXu^{-1} \equiv un(a_1,\ldots,a_k)u^{-1} \equiv n(b_1,\ldots,b_p) \equiv X \mod \nl_{\geq p+1}.\] 
 Thus by (\ref{stingev2}) 
 \[ un(a_1,\ldots,a_k)u^{-1} \equiv [X,n_p]+n(a_1+d_1x^p,\ldots,a_k+d_kx^p) \mod \nl_{\geq p+2},\]
 with $d_i := d_{\alpha_i}$ and $n_p\in \nl_p$. Thus \[n(a_1+d_1x^p,\ldots,a_k+d_kx^p)-X \equiv n(b_1,\ldots,b_k)-X \mod [X,\nl_p].\] Since $f: \nl_{p+1}\rightarrow \F$ is a linear map with kernel $[X,\nl_p]$,
 \[ f(b_1-a_1,\ldots,b_k-a_k) = f(d_1x^p,\ldots,d_kx^p)= x^pf(d_1,\ldots,d_k).\]
 Define $c := f(d_1,\ldots,d_k)$. Then $f(b_1-a_1,\ldots,b_k-a_k)\in c\F^{(p)}$ if and only if $n(a_1,\ldots,a_k)$ is $U$-conjugated modulo $\nl_{\geq p+2}$ to $n(b_1,\ldots,b_k)$.\\
 Since being $U$-conjugated modulo $\nl_{\geq p+2}$ is an equivalence relation, we have that $c\in \F^{(p)}$. Now $c\not=0$, because by Proposition \ref{Carpros} over an algebraically closed field the orbit of $X$ in $\nl$ contains $X+\bigoplus_{i=2}^{ht(R)} \nl_i$. Thus $c\in (\F^\times)^p$.
\end{proof}

\begin{lem}\label{ssadc}
 If $G$ is simple of adjoint type and $p$ is bad for $G$, then there exists a bad pair for $G$.
\end{lem}
\begin{proof}
Choose an $\alpha\in \Delta$, define \[U_{\hat{\alpha}} := \prod_{\gamma\in\{\gamma\in R^+-\{\alpha\}\mid ht(\gamma)\leq p-1\}} U_\alpha.\]
Since $G$ is of adjoint type, $\Phi : X_*(T)\rightarrow \hom_\Z(\Z R(G,T),\Z)$ is surjective. Thus for every $n\in \nl'$ there is exactly one $t\in T$ such that $tnt^{-1} \in X+\nl_{\geq2}$. Therefore by the proof of Proposition \ref{linfun}, for every $n\in \nl'$, there exists an unique $b\in TU_{\hat{\alpha}}$ such that $bnb^{-1} = X+n(a_1,\ldots,a_k)+n_{p+2}$, with $n_{p+2}\in \nl_{\geq p+2}$. Write $n = \sum_{\alpha\in R} x_\alpha E_\alpha$. The $a_1,\ldots,a_k$ depend polynomially on $x_\alpha$ for $\alpha\in R_i^+$ with $i\leq p+1$. Let $f_i$ be the polynomials such that $a_i = f_i(x_\alpha)$. The $f_i$ are homogeneous of degree $-p$:
\begin{align*}
bnb^{-1} &\equiv X+n(a_1,\ldots,a_k) &&\mod \nl_{\geq p+2}\\
b\lambda n b^{-1} &\equiv \lambda (X+n(a_1,\ldots,a_k)) &&\mod \nl_{\geq p+2}\\
tb\lambda n b^{-1} t^{-1} &\equiv X+\frac{\lambda}{\lambda^{p+1}}n(a_1,\ldots,a_k) &&\mod \nl_{\geq p+2},
\end{align*}
where $t\in T$ is such that $\gamma(t)=\frac{1}{\lambda}$ for all $\gamma\in \Delta$.\\
Define $\chi(n) := f(f_1(x_\alpha),\ldots,f_k(x_\alpha))$ for $n\in \nl$.\\
Choose a $g:\F\rightarrow \F^k$ to be a right inverse of $f$, ie. $fg=id$, such that $n(g(\co))\subset \nl_{p+1}(\co)$. Define $\eta : \F\rightarrow \nl'$ by $\eta(a) := X+n(g(a))$. Now $(\chi,\eta)$ is a bad pair for $G$.
\end{proof}

\begin{lem}\label{adgng}
 If $(\eta,\chi)$ is a bad pair for $Ad(G)$, then there exists a bad pair for $G$.
\end{lem}
\begin{proof}
 Let $Ad : G\rightarrow Ad(G)$ be the natural morphism. By Lemma \ref{darnib},  $d(Ad):\nl\rightarrow \nl^{ad}$ is a bijection, let $da : \nl^{ad}\rightarrow \nl$ be its inverse. If $n,n'\in \nl$ are conjugated by $G$, then their image is conjugated by $Ad(G)$. If $E_\alpha$ is a Chevalley basis for $G$, then $d(Ad)(E_\alpha)$ is one for $Ad(G)$. Thus $(da\circ \eta, \chi\circ d(Ad))$ is a bad pair for $G$.
\end{proof}

\begin{ste}\label{spottedbp}
Assume that $G$ is a $\F$-split reductive group and $p$ is bad for $G$. Then there exists a bad pair for $G$.
\end{ste}
\begin{proof}
By Lemma \ref{adgng} we may assume that $G$ is semi-simple of adjoint type. Assume that $G = G_1\cdots G_m$, with $G_i$  the simple connected normal subgroups of $G$ and $p$ bad for $G_1$. Let $(\eta_1,\chi_1)$ be the bad pair of Lemma \ref{ssadc}. Define $\chi(n_1+\ldots+n_m) := \chi_1(n_1)$, for $n_i\in\nl^i:= \nl\cap \gl_i$. The function $\chi$ is well-defined, since $\nl = \nl^1\oplus \cdots\oplus \nl^m$. \\
Assume that $n,n'\in\nl'$ are conjugated, then there exists a $b\in B$ such that $bnb^{-1} = n'$. Write $b= t \prod_{i=1}^m u_i$ with $u_i\in U_i= U\cap G_i$ and $t\in T$. Then $u_1tn_1(u_1t)^{-1}=n'_1$. Since $G=G_1\times \cdots\times G_m$, $n_1$ and $n'_1$ are also conjugated by $G_1$. We conclude that if $n$ is conjugated with $n'$, then $\chi(n)\equiv \chi(n') \mod \F^{(p)}$. For $2\leq i\leq m$, let $n_i\in\nl^i$ be a regular nilpotent element. Define $\eta(x) := \eta_1(x)+n_2+\cdots+n_m$. Thus $(\eta,\chi)$ is a bad pair for $G$.
\end{proof}

\begin{ste}\label{bnhc}
 If $G$ is a $\F$-split reductive group and char $\F$ is bad, then there are infinitely many nilpotent orbits and Howe's conjecture on the Lie algebra does not hold.
\end{ste}
\begin{proof}
 This follows from Theorem \ref{spottedbp} and Theorem \ref{bpnhc}.
\end{proof}

\subsubsection{The example $SO_5(\F)$, char $\F=2$}
In this section $\F$ has characteristic $2$. We follow \cite[\S 7.4.7(6)]{Sp98} for the definition of $SO_5(\F)$. Let $V = \F^5$ and let $Q$ be the quadratic form on $V$ defined by
\[Q(e_0,e_1,e_2,e_3,e_4) := e_0^2+e_1e_3+e_2e_4.\]
Now we define $SO_5(\F)$ to be the subgroup of $t\in GL(V)$ with $Q(tv)=Q(v)$ for all $v\in V$. Then
\[ T := \{ t(t_1,t_2):=\left(\begin{array}{ccccc} 1&0&0&0&0\\0&t_1&0&0&0\\0&0&t_2&0&0\\0&0&0&t_1^{-1}&0\\0&0&0&0&t_2^{-1}\end{array}\right) : t_i\in\F^\times \} \]
is a maximal torus of $SO_5$ that is $\F$-split. Define, for $i=1,2$, the character $\eps_i$ of $T$ by
\[ 
\eps_i(t(t_1,t_2)):= t_i.
\]
Then $R(G,T) = \{ \pm \eps_i, \pm \eps_i\pm\eps_j \mid i\not=j\}$. Let $R^+ := \{\eps_1-\eps_2,\eps_2,\eps_1,\eps_1+\eps_2\}$ be a system of positive roots and $\Delta := \{ \eps_1-\eps_2,\eps_2\}$ the corresponding set of simple roots. We get the following positive root spaces in the Lie algebra $\gl$ of $SO_5$:
\begin{align*}
 \gl_{\eps_1-\eps_2} &:= \left(\begin{array}{ccccc} 0&0&0&0&0\\0&0&c&0&0\\0&0&0&0&0\\0&0&0&0&0\\0&0&0&c&0\end{array}\right)
 &\gl_{\eps_1} &:= \left(\begin{array}{ccccc} 0&0&0&a&0\\0&0&0&0&0\\0&0&0&0&0\\0&0&0&0&0\\0&0&0&0&0\end{array}\right)\\
 \gl_{\eps_2} &:= \left(\begin{array}{ccccc} 0&0&0&0&b\\0&0&0&0&0\\0&0&0&0&0\\0&0&0&0&0\\0&0&0&0&0\end{array}\right)
 &\gl_{\eps_1+\eps_2} &:= \left(\begin{array}{ccccc} 0&0&0&0&0\\0&0&0&0&e\\0&0&0&e&0\\0&0&0&0&0\\0&0&0&0&0\end{array}\right)
\end{align*}
Thus $X := \left(\begin{array}{ccccc} 0&0&0&0&1\\0&0&1&0&0\\0&0&0&0&0\\0&0&0&0&0\\0&0&0&1&0\end{array}\right)$. Also $\nl_2 = \gl_{\eps_1}$ and $\nl_3 = \gl_{\eps_1+\eps_2}$. The function $[X,\cdot] :\nl_2\rightarrow \nl_3$ is as follows:
\[
\left[ \left(\begin{array}{ccccc} 0&0&0&0&1\\0&0&1&0&0\\0&0&0&0&0\\0&0&0&0&0\\0&0&0&1&0\end{array}\right) ,  \left(\begin{array}{ccccc} 0&0&0&a&0\\0&0&0&0&0\\0&0&0&0&0\\0&0&0&0&0\\0&0&0&0&0\end{array}\right) \right] = 0.
\]
Thus according to Proposition \ref{linfun} and its proof, $X+eE_{\eps_1+\eps_2}$ is $U$-conjugated with $X+e'E_{\eps_1+\eps_2}$ if and only if $e\equiv e' \mod \F^{(2)}$. Now we follow Lemma \ref{ssadc}. We take $U_{\hat{\eps_2}} := U_{\eps_1-\eps_2}$. Define for $a,b,c,d\in \F$
\[ n(a,b,c,e) := \left(\begin{array}{ccccc} 0&0&0&a&b\\0&0&c&0&e\\0&0&0&e&0\\0&0&0&0&0\\0&0&0&c&0
                       \end{array}\right). \]
Assume that $b,c\not=0$, then there is by the Lemma a unique $g\in TU_{\hat{\eps_2}}$ such that $gn(a,b,c,e)g^{-1} = X + e'E_{\eps_1+\eps_2}$ for some $e'\in \F$. Lets compute $e'$: first get the $b$ and $c$ to $1$ by conjugating with $t := t(c^{-1}b^{-1}, b^{-1})$, then
\[tn(a,b,c,e)t^{-1} = n(ab^{-1}c^{-1},1,1,ec^{-1}b^{-2}).\]
By conjugating the result with $\left(\begin{array}{ccccc} 1&0&0&0&0\\0&1&ab^{-1}c^{-1}&0&0\\0&0&1&0&0\\0&0&0&1&0\\0&0&0&ab^{-1}c^{-1}&1\end{array}\right)$ we get $n(0,1,1,ec^{-1}b^{-2})$. Thus $e' = \frac{e}{bc^{-1}}$. Assume that $b',c'\not=0$. Thus $n(a,b,c,e)$ is conjugated with $n(a',b',c',e')$ if and only if $\frac{e}{cb^2} \equiv \frac{e'}{c'b'^2}\mod \F^{(2)}$.

\section{Howe's conjecture and $\kappa_v(G)$}
In this section we assume that $p$ divides $\kappa_v(G)$, ie the characteristic of $\F$ divides the cokernel of the map:
\begin{align*}
 \Phi : X_*(T)&\rightarrow \text{Hom}_\Z(\Z R(G,T),\Z)\\
 \gamma &\mapsto (\alpha \mapsto \left<\gamma,\alpha\right>)
\end{align*}

We will follow the same strategy as in section \ref{redbp}.
By the proof of Proposition \ref{ubaan} there exists integers $z_i\in\Z$ such that $\chi: \nl' \rightarrow \F^\times$ defined by $\chi(n):= \prod_{i=1}^m n_{\alpha_i}^{z_i}$ is surjective and $\chi :\nl'\rightarrow \F^\times/(\F^\times)^p$ is $B$-invariant. Take a $\eta : \F^\times \rightarrow \nl$ such that $\eta$ is algebraic and $\chi\eta$ is the identity. By the proof of Proposition \ref{ubaan} we can choose $\eta$ in such a way that for all $\alpha\in \co^\times$: $\eta(\alpha)_{\gamma}\in \co^\times $ for all $\gamma\in \Delta$ and $\eta(\alpha)_{\gamma}=0$ for all $\gamma \in R-\Delta$. The functions $\chi,\eta$ play the role of bad pair in this case.

\begin{lem}\label{kiesN}
 There exists a $N>0$ such that for all $n\in\N_{>0}$, $k\in K$ and $\alpha\in \co^\times$:
\[k\eta(\alpha)k^{-1} \in \nl + \pi^{Nn}L \Rightarrow k\in (B\cap K)K_n.\]
\end{lem}
\begin{proof}
 See the proof of Lemma \ref{kiesN3}.
\end{proof}
Define the following $B$-invariant open set of $\nl$:
\[ V_{\alpha,s} := \left\{n\in \nl' \mid \chi(n)\equiv \alpha \mod (1+\pi^s\co)(\F^\times)^p \right\}. \]

Define $n(a_1,\ldots,a_m) := \sum_{i=1}^m a_iE_{\alpha_i}$.\\
Define $z := \sum_{i=1}^m z_i$. Then $\chi(\pi^n x) = \pi^{zn}\chi(x)$, for all $x\in \gl$ and $n\in \Z$.
\begin{lem}\label{aebmss}
Let $n\in\N_{>0}$, $\alpha\in\co^\times$ and $\beta\in \co^\times$.\\
If
 \[ \int_{V_{\pi^{-znN}\beta,n}}\int_{k\in K} 1_{\pi^{-Nn}\eta(\alpha)+L}(k(X+Z)k^{-1})dkdZdX>0,\]
then $\alpha \equiv \beta \mod (1+\pi^n\co)(\co^\times)^p$.
\end{lem}
\begin{proof}
So there exist a $k\in K$ and $l\in L$ such that \[k\pi^{-nN}\eta(\alpha)k^{-1}+l\in V_{\pi^{-znN}\beta,n}+\nl_{\geq2}\subset \nl.\] Thus $k\in (K\cap B)K_n$ by Lemma \ref{kiesN}, because $k\eta(\alpha)k^{-1}\in \nl+\pi^{nN}L$. Take $b_k\in K\cap B$ and $k_n\in K_n$ such that $k=k_nb_k$. Take $a_1,\ldots,a_m\in \co^\times$ and $n_2\in \nl_{\geq 2}(\co)$ such that $b_k\eta(\alpha)b_k^{-1}=n(a_1,\ldots,a_m)+n_2$. By the construction of $\chi$, there exists a $\gamma\in\F^\times$ such that $\chi(n(a_1,\ldots,a_m)) = \alpha\gamma^p$.  Since $k_n\in K_n$ and $n(a_1,\ldots,a_m)+n_2\in L$, there exists a $l'\in L$ such that\\ $k_n(n(a_1,\ldots,a_m)+n_2)k_n^{-1}= n(a_1,\ldots,a_m)+n_2+\pi^nl'$. Thus
\begin{align*}
 \chi(k\eta(\alpha) k^{-1}+\pi^{nN}l) &=\chi(k_nb_k\eta(\alpha)b_k^{-1}k_n^{-1}+\pi^{nN}l)\\
 &= \chi(n(a_1,\ldots,a_m)+n_2+\pi^n l'+\pi^{nN}l) = \prod_{i=1}^m (a_i+\pi^nl_i)^{z_i},
\intertext{for some $l_i\in \co$. Since the $a_i$ are in $\co^\times$,}
 &\prod_{i=1}^m (a_i+\pi^nl_i)^{z_i}\equiv \prod_{i=1}^m a_i^{z_i}  = \alpha\gamma^p\mod (1+\pi^n\co).
\intertext{Because $\chi(\pi^{-nN}x)=\pi^{-znN}\chi(x)$ for all $x\in \gl$,}
 \chi(k\pi^{-nN}\eta(\alpha)k^{-1}+l) &\equiv (\alpha\gamma^p)\pi^{-znN} \mod (1+\pi^{n}\co). 
\intertext{
Since $k\pi^{-nN}\eta(\alpha)k^{-1}+l\in V_{\pi^{-znN}\beta,n}+\nl_{\geq 2}$,} \chi(k\pi^{nN}\eta(\alpha)k^{-1}+l) &\equiv \pi^{-znN}\beta \mod (\F^\times)^p(1+\pi^{n}\co).\\ 
\intertext{Thus}
 \pi^{-znN}\beta &\equiv \chi(k\pi^{-nN}\eta(\alpha)k^{-1}+l)\equiv \pi^{-znN}\alpha \mod ((\F^\times)^p (1+\pi^n\co).
\end{align*}
Then $\alpha \equiv \beta \mod (\F^\times)^p(1+\pi^n\co)$. Because $(\F^\times)^p\cap \co^\times = (\co^\times)^p$ and $\alpha,\beta\in \co^\times$, the Lemma follows.
\end{proof}

\begin{ste}\label{scnhc}
 Let $G$ be a $\F$-split reductive group. Assume that $\text{char}\;\F | \kappa_v(G)$, then Howe's conjecture does not hold.
\end{ste}
\begin{proof}
 The proof is simular to the one of Theorem \ref{bpnhc}.\\
 Let $\alpha_1,\ldots,\alpha_k$ be representatives of the cosets of $(1+\pi^n\co)\co^\times$ in $\co^\times$. Define for $1\leq i\leq k$ the following distribution and function:
 \begin{align*}
  D_i(f) &:= D_{V_{\pi^{-znN\alpha_i,n}}}(f)\\
  f_i &:= 1_{\pi^{-nN\eta(\alpha_i)}+L}
 \end{align*}
 Let $c_i := D_i(f_i)>0$. By Lemma \ref{aebmss} $D_i(f_j) = c_i\delta_{ij}$. Therefore $\dim J_L(\omega)\geq k$.
\end{proof}

\section{Howe's conjecture in good characteristic}
Howe's conjecture does not hold when the characteristic is bad or $p|\kappa_v(G)$. In this section we investigate Howe's conjecture in good characteristic. Throughout this section we assume that $p$ is good for $G$.

\subsection{Associated cocharacters to nilpotent elements}
In this subsection we recall the theory of associated cocharacters.
Let $\tau \in X_*(G)$. For $z\in \Z$, we define the following subspaces of $\gl$:
\begin{align*}
\gl(z;\tau) &:= \{ X\in \gl \mid \forall[a\in \F]\;\tau(a)X\tau(a)^{-1} = a^zX\}\\
\gl(\geq z;\tau) &:= \bigoplus_{i\geq z} \gl(i;\tau). 
\end{align*}
We sometimes abbreviate $\gl(z;\tau)$ ($\gl(\geq z;\tau)$) by $\gl(z)$ ($\gl(\geq z)$ resp.), in which case the cocharacter $\tau$ should be clear from the context.\\
A nilpotent element $X\in\gl$ is called distinguished if each torus contained in $Z_G(X)$ is contained in the center of $G$.\\
A cocharacter $\tau$ of $G$ is called associated to $X$ if $X \in \gl(2,\tau)$ and if there exists a Levi subgroup $L$ in $G$ such that $X$ is distinguished nilpotent in $\mathfrak{l}$ and such that $\text{im }\tau \subset (L,L)$.\\
Let $X\in\gl$ be nilpotent, define $N(X) := \{ g\in G \mid Ad(g)X\in \F X\}$.
\begin{lem}\cite[Lemma 25]{Mc04}
 Let $S$ be any maximal torus of $N$. Then there is a unique cocharacter in $X_*(S)$ associated with $X$.
\end{lem}
\begin{ste}\cite[Theorem 26]{Mc04}
 Let $X\in \gl$ be nilpotent. Assume that the $G$-orbit of $X$ is separable. Then there exists a cocharacter $\tau$ associated to $X$ which is defined over $\F$.
\end{ste}
Let $\tau$ be a cocharacter associated to $X$, we define
\begin{align*}
\pl_X &:= \gl(\geq 0;\tau)\\
\nl_X &:= \gl(\geq 1;\tau)
\end{align*}
The Lie algebras $\pl_X$ and $\nl_X$ are independent of the choice of $\tau$.
\subsection{First proof of Howe's conjecture}
\begin{lem}\label{satopr59c}
 Suppose that char$(\F)$ is good for $G$ and $\F$ is algebraically closed. Let $X$ be nilpotent. Let $\lambda$ be a cocharacter associated with $X$. Then
 \[ [\gl(-1),X] = \gl(1)\]
 and
 \[[\nl_X,X] = \gl(\geq 3).\]
\end{lem}
\begin{proof}
We follow the same line as the proof of \cite[Proposition 5.9(c)]{Ja04}.\\
Let $G$ be a group satisfying the standard hypotheses:
\begin{enumerate*}
 \item The derived group of $G$ is simply connected.
 \item The characteristic of $\F$ is good for $G$.
 \item There exists a $G$-invariant nondegenerate bilinear form on $\gl$.
\end{enumerate*}

By \cite[Proposition 5.8 and Lemma 5.7]{Ja04}
\[[\gl(-1),X]=\gl(1)\]
and
\[[\nl_X,X] = \gl(\geq 3).\]

Now we show that the Lemma holds for $G$ if and only if it holds for $G_{\text{der}}$.\\
The cocharacter $\tau$ associated to $X$ in $G$ is also the cocharacter $\tau$ associated to $X$ in $G_{\text{der}}$. Also $\gl(-1),\gl(1),\gl(\geq 3)\subset \gl'$.\\

When $G$ is simply connected and the characteristic is very good, then $G$ satisfies the standard hypotheses. The Lemma holds for $GL_n$ by \cite[Lemma 2]{Ho74}, thus for $SL_n$ as well. Therefore the Lemma holds for all simply connected groups in good characteristic. Hence also for products of those groups.\\

Let $G = R(G)G_1,\cdots,G_m$ with $G_i$ the simple normal connected subgroups of $G$. Let $G_i'$ be the simply connected group belonging to $G_i$. Let $\pi : R(G)\prod_{i=1}^m G_i' \rightarrow G$ be the natural surjective homomorphism. Now $d\pi$ is surjective on the nilpotent elements and maps the associated cocharacter of a nilpotent element to the associated cocharacter of its image. Since the Lemma holds for $R(G)\prod_{i=1}^m G_i'$ it also holds for $G$.
\end{proof}

\begin{ste}\label{Higc}
 Let $G$ be a reductive group and the characteristic of $\F$ be good for $G$. If the nilpotent orbits of $G$ in $\gl$ are separable, then Howe's conjecture holds.
\end{ste}
\begin{proof}
Basically the proof of Harish-Chandra in \cite{HC99} does the job. The proof of Harish-Chandra is in the characteristic zero case. We will only mention the two adjustments to make it work in this case as well. The adjustments are all in the proof of \cite[Theorem 13.1]{HC99}.\\
Let $X_0$ be a nilpotent element of $\gl$.\\
In \cite[\S13.1]{HC99} Harish-Chandra completes $X_0$ to a Jacobson-Morosow triple. The analogue in positive characteristic is of course the cocharacter associated to $X_0$. By \cite[Theorem 26]{Mc04} there exists a cocharacter $\phi$ associated to $X_0$ which is defined over $\F$, because the orbit of $X_0$ is separable. In the proof of \cite[Proposition 34]{Mc04}:\\
 "Well, by \cite[Proposition 5.9(c)]{Ja04}, we have $\overline{\text{Ad}(P)X}=\oplus_{i\geq 2} \gl(i,\phi)$. Since the orbit of $X$ is separable, the differential of the orbit map is surjective."\\
Thus $\gl(\geq 2)\subset [\gl,X]$. Therefore with Lemma \ref{satopr59c} we have $\nl_X \subset [X,\gl]$. (This is Lemma 13.2 of \cite{HC99}).\\
The proof of Lemma 13.5 of \cite{HC99} uses the exponential map. We can replace the exponential map by the mock exponential map of Adler \cite{Ad98}. See in particular \cite[Proposition 1.6.3]{Ad98}.
\end{proof}

Observe that the conditions of Theorem \ref{Higc} are geometric conditions; they only depend on the algebraic group and the algebraic closure of $\F$, not on the $\F$-form of $G$.

\begin{gev}\label{vghc}
 If $G$ is a simple group and $\F$ is very good for $G$, then Howe's conjecture holds.
\end{gev}
\begin{proof}
 By \cite{Ri67} the nilpotent orbits are separable.
\end{proof}

\subsection{The case $SO_{3}(\F)$ (char $\F=2$)}
In this section char $\F=2$.\\

Although there are infinitely many nilpotent conjugacy classes in $SO_3(\F)$ and the nilpotent orbits are not separable, Howe's conjecture holds for $SO_3(\F)$. We again follow \cite{HC99}, but have to make a few more modifications.\\

The next lemma and its proof are \cite[Lemma 12.2]{HC99}, with $\nl^G$ instead of $\cn$.
\begin{lem}\label{L12.2}
 Let $\omega \subset \mathfrak{g}$ be a compact set.\\
 Let $S$ be a split torus and $K$ the stabilizer of $0$ in the apartment of $S$ (in the extended building). Take $\Phi^+$ a system of positive roots of $(G,S)$. Let $\nl$ be the Lie algebra for $U^+$, $\onl$ be the Lie algebra for $U^-$ and $\ml$ the Lie algebra of $M := Z_G(S)$.\\
 There is a lattice $\Lambda$ such that
 $$Ad(G)\omega = \Lambda + Ad(KS)(\nl\cap \Lambda)$$
\end{lem}
\begin{proof}
By Bruhat-Tits one has
$$G = KSFK$$
for some finite subgroup $F$ of $M$.\\
Since $\gl = \onl\oplus\ml\oplus\nl$ one has compact subsets $\omega_1,\omega_2,\omega_3$ in $\onl$, $\ml$ and $\nl$ respectively, such that 
$$Ad(FK)\omega\subset \omega_1\oplus \omega_2\oplus \omega_3.$$
Hence $Ad(G)\omega \subset Ad(KS)(\omega_1\oplus \omega_2\oplus \omega_3)$.\\
Now $Ad(S)\omega_1$ is contained in a compact lattice of $\onl$, since $v(\alpha(s))\geq 0$ for all $\alpha\in\Phi^-$ and $Ad(S) \omega_2=\omega_2$. Therefore there is a lattice $L$ such that
\[Ad(G)\omega \subset Ad(K) \left(L + Ad(S)(\nl\cap L)\right).\qedhere\]
\end{proof}
Since $^G\nl=\cn$ in characteristic $0$, Lemma 12.2 of Harish-Chandra works with $\cn$. For the group $SO_3(\F)$ this is not the case. Therefore we shall work with $^G\nl$ instead of $\cn$. We start with the definition of $SO_3(\F)$.\\
Define $Q(e_0,e_1,e_2) := e_0^2+e_1e_2$.

\[ SO_3(\F) := \{ g\in GL_3 \mid Q(gv)=Q(v)\}\]
Let $\gamma$ be the following cocharacter of $SO_3$.
$$\gamma(t) := \left(\begin{array}{ccc} 1&0&0\\0&t&0\\0&0&t^{-1}\end{array}\right)$$
Let $T$ be the following subgroup of $SO_3$:
\[ T := \{ \gamma(t) : t\in \F^\times\}\]
Now $T$ is a maximal torus of $SO_3$.\\
The Lie algebra of $SO_3$ is of the following form:
\[\gl := \{ \left(\begin{array}{ccc} 0&a&b\\0&c&0\\0&0&c\end{array}\right) : a,b,c\in\F\}. \]
With respect to the cocharacter $\gamma$ we have a decomposition of the Lie algebra:\\ $\gl := \gl(-1)\oplus \gl(0) \oplus \gl(1)$ with
\begin{align*}
 \nl&:= \gl(1) = \{\left(\begin{array}{ccc} 0&0&b\\0&0&0\\0&0&0\end{array}\right) : b\in \F\}\\
 \tl &:= \gl(0) = \{\left(\begin{array}{ccc} 0&0&0\\0&c&0\\0&0&c\end{array}\right) : c\in \F\}\\
 \onl &:= \gl(-1) =  \{\left(\begin{array}{ccc} 0&a&0\\0&0&0\\0&0&0\end{array}\right) : a\in \F\}
\end{align*}
Take on $\gl$ the following norm:
$$ \left|\left(\begin{array}{ccc} 0&a&b\\0&c&0\\0&0&c\end{array}\right)\right|= \max(|a|,|b|,|c|)$$

For the extended version of Howe's conjecture, Harish-Chandra needs to consider all nilpotent orbits. But for the regular Howe's conjecture we can restrict ourselves to one nilpotent orbit, namely the orbit of $$n := \left(\begin{array}{ccc}0&0&1\\0&0&0\\0&0&0\end{array}\right).$$

Let $N := {^Gn}\cup \{0\} = {^G\nl}$.\\
Define, for $a,b\in\F$, the following elements of $\mathfrak{so}_3$ and $SO_3$:
\[
n_{a,b} := \left(\begin{array}{ccc} 0&a&b\\0&0&0\\0&0&0\end{array}\right),\;
u_b := \left(\begin{array}{ccc} 1&0&b\\0&1&b^2\\0&0&1\end{array}\right),\;
\omega := \left(\begin{array}{ccc} 1&0&0\\0&0&1\\0&1&0\end{array}\right).
\]

\begin{lem}\label{clasN}
 $N = \{ n_{a,b} \mid \exists(y\in \F)\; y^2=ab\}$
\end{lem}
\begin{proof}
The conjugation action of the generators of $SO_3(\F)$ on the nilpotent elements is as follows:
\begin{align*}
 u_cn_{a,b}u_c &= n_{a,c^2a+b}\\
 \gamma(x)n_{a,b}\gamma(x)^{-1} &= n_{x^{-1}a,xb}\\
 \omega n_{a,b} \omega &= n_{b,a}
\end{align*}
The lemma follows after some calculations.
\end{proof}
\begin{gev}\label{Ncc}
 The set $N$ is closed in $\gl$ and $cN = N$ for all $c\in \F^\times$.
\end{gev}
\begin{proof}
 The nilpotent elements are closed in $\gl$. The function $Q :n_{a,b}\mapsto ab$ is a continuous function form $\mathcal{N}$ to $\F$. Since $\F^{(2)}$ is closed in $\F$, so it $Q^{-1}(\F^{(2)})$. The latter is equal to $N$ by Lemma \ref{clasN}. Since closed sets of closed subspaces are closed, $N$ is closed in $\gl$. The second statement is obvious.
\end{proof}
\begin{lem}\label{hggvN}
 Let $X \in N$.\\
 There is a cocharacter $\tau$ such that $X\in \gl(1)$ and $\gl(1)\subset [X,\gl]$.
\end{lem}
\begin{proof}
 Since these statements are $G$-invariant, we may and will assume that $X= n_{0,1}$. In this case take $\tau := \gamma$. Clearly $X\in\nl\subset \gl(1)$. Now
 $$\left(\begin{array}{ccc} 0&0&0\\0&c&0\\0&0&c\end{array}\right)\left(\begin{array}{ccc} 0&0&1\\0&0&0\\0&0&0\end{array}\right) +\left(\begin{array}{ccc} 0&0&1\\0&0&0\\0&0&0\end{array}\right)\left(\begin{array}{ccc} 0&0&0\\0&c&0\\0&0&c\end{array}\right)= \left(\begin{array}{ccc} 0&0&c\\0&0&0\\0&0&0\end{array}\right),$$
thus $\nl \subset [X,\gl]$.
\end{proof}
\begin{ste}
 Howe's conjecture holds in $SO_3(\F)$.
\end{ste}
\begin{proof}
 We follow Harish-Chandra \cite{HC99} again and mention the adjustments. We replace $\cn$ by $N=\;^G\nl$. The proof of Harish-Chandra uses three properties of $\cn$ (in brackets the Lemma's in \cite{HC99} where the property is used):
 \begin{enumerate*}
  \item $\cn\cap S$ is compact (Lemma 11.9)
  \item For all compact subsets $\omega$ in $\gl$ there exists a lattice $L_1$, such that $^G\omega \subset L_1+\cn$. (Lemma 12.2)
  \item If $c\in \F$ and $Y\in \cn$, then $cY\in \cn$ (Lemma 12.3)
 \end{enumerate*}
 By Corollary \ref{Ncc} (1) and (3) also hold for $N$ and (2) is Lemma \ref{L12.2}. Now we are left to prove Theorem 13.1 for $X\in N\cap S$. By Lemma \ref{hggvN} we can use the proof of Theorem \ref{Higc}.
\end{proof}
This example shows that the separability of the nilpotent orbits is not a necessary condition for Howe's conjecture to hold.

\subsection{The case $PGL_n (\F)$ with char $\F|n$}
In this section we generalize the results in the previous section to the group $PGL_n\;(\F)$. This is the group of $\F$-points of the algebraic quotient of $GL_n$ by its centrum of diagonal matrices $Z$. We have the exceptional isomorphism $PGL_2 \cong SO_3$. Let $G := PGL_n$. We identify $\gl$ with $\Gl_n/\zl$. Now $\Gl_n/\zl := \{ X+\zl : X\in \Gl_n\}$. Define $p := \text{char}\;\F$. The nilpotent elements of $\gl$ are exactly those $X+\zl$ such that $X^{p^n} \in \zl$. We define the following $G$-invariant function $\phi$ on $\mathcal{N}$: for $X\in \gl(\F)$ let $a\in \F$ be such that $X^{p^n} = aI_n$, with $I_n$ the identity matrix. Then $\phi(X+\zl) := a+\F^{(p^n)}$. If $X+\zl=X'+\zl$, then $X-X'\in \zl(\F)$. Thus $\phi$ is well defined.
\begin{lem}
The following statements hold for $\phi$:
\begin{enumerate*}
 \item $\phi$ is $G$-invariant
 \item $\F^{(p^{n-1})} \subset \text{Im}\; \phi$
 \item Let $X+\zl$ be a nilpotent element of $\gl$. Then $\phi(X) \in \F^{(p^n)}$ if and only if there exists a nilpotent matrix $n\in\Gl_n$ such that $n\in X+\zl$.
\end{enumerate*}
\end{lem}
\begin{proof}
1. trivial.\\
2. Let $M_x$ be a block-diagonal matrix consisting of $\frac{n}{p}$ blocks with on each $(p\times p)$-block the matrix
\[
\left(\begin{array}{ccccc}
0 & 0& \cdots & 0&x\\
1&0&0&\cdots &0\\
0&1&0&\cdots&\vdots\\
0&0&\ddots &0&0\\
0&\cdots&0&1&0
\end{array}\right)
\]
Then $M_x^p= xI_n$ thus $M_x^{p^n} = x^{p^{n-1}}I_n$.\\
3. For a nilpotent matrix $n\in\Gl_n$ we have that $\phi(n+\zl)=0$, thus the only if part is clear. Assume that $X^{p^n} = a^{p^n}I_n$, then $(X-aI_n)^{p_n} = X^{p^n}-a^{p^n}I_n=0$. Thus $X-aI_n$ is nilpotent.
\end{proof}
\begin{gev}\label{pglio}
 The number of nilpotent orbits is infinite.
\end{gev}
\begin{proof}
 The group $\F^{(p^{n-1})}/\F^{(p^n)}$ is as group isomorphic to $\F/\F^{(p)}$, $\F/\F^{(p)}$ is infinite and $\F^{(p^{n-1})} \subset \text{Im}\; \phi$.
\end{proof}
Thus not all nilpotent orbits are separable. In fact the orbit of \[ x := \left(\begin{array}{cccc} 0&1&0&0\\ 0&0&\ddots&0\\0&0&0&1\\0&\cdots & &0\end{array}\right),\]
the superdiagonal entries of $x$ are $1$, is not separable,
since the commutator with
\[ x' := \left(\begin{array}{cccc} 0&0&0&0\\ 1&0&0&0\\0&\ddots&0&0\\0&0&n&0\end{array}\right),\]
the subdiagonal entries of $x'$ are from left to right equal to $1,2,\cdots,n$, is equal to (char $\F|n$)
\[
 \left(\begin{array}{cccc} 1&0&0&0\\0&\ddots &0&0\\0&0&1&0\\0&0&0&-(n-1)\end{array}\right) = Id_n\in\zl.
\]
\begin{lem}\label{ldhpgln}
 If $p|n$, then $\{H_\alpha : \alpha\in \Delta\}$ are linearly dependent. 
\end{lem}
\begin{proof}
 Let $T$ be the torus of diagonal matrices. For $i=1,\ldots,n$, define
 \[ \eps_i\left(\begin{array}{ccc} x_1&0&\cdots \\0&\ddots&0\\ \cdots&0&x_n\end{array}\right) = x_i.\]
 Let $\Delta = \{ \eps_1-\eps_2,\ldots,\eps_{n-1}-\eps_n\}$, then
 \[\sum_{i=1}^{n-1} H_{\eps_i-\eps_{i+1}} = [x,x']=0.\qedhere\]
\end{proof}

Define $N := \{ x+\zl : x\in\Gl_n \mid x \text{ is nilpotent}\}$.

\begin{gev}\label{NccN}
 $N$ is a closed subset of $\gl$ and $cN=N$ for $c\in \F^\times$.
\end{gev}
\begin{proof}
The map $\phi$ is continuous and $0\in \F/\F^{(p^n)}$ is closed. Thus $N = \phi^{-1}(0)$ is closed in $\mathcal{N}$. Because $\mathcal{N}$ is closed in $\gl$, so is $N$.\\
If $x\in \gl$ is nilpotent, then, for all $c\in \F$, $cx$ is also nilpotent. Thus $cN=N$.
\end{proof}

\begin{lem}\label{ecax}
 For every nilpotent element $X\in \Gl_n$ there exists a cocharacter $\gamma$, such that $X\in \nl(\gamma)$ and $\nl(\gamma) \subset [X,\mathfrak{p}(\gamma)]$.
\end{lem}
\begin{proof}
 We follow \cite{Ho74} and its notation. See page 311 of loc. cit.. Define for $x\in \F$ the element $\gamma(x)\in M$ to be the transformation which acts on $C_i$ by multiplication by $x^i$. Then $\mathcal{U}= \nl(\gamma)$ and by \cite[Lemma 2]{Ho74} $\nl(\gamma) \subset [X,\mathfrak{p}(\gamma)]$.
\end{proof}
\begin{gev}\label{gecax}
 For every nilpotent element in $N$ there exists a cocharacter $\gamma$, such that $X\in \nl(\gamma)$ and $\nl(\gamma) \subset [X,\mathfrak{p}(\gamma)]$.
\end{gev}
\begin{proof}
 Let $X\in GL_n$ and let $\gamma \in X_*(G)$ be the cocharacter of Lemma \ref{ecax}. Let $\varphi : GL_n\rightarrow PGL_n$ be the natural homomorphism. Because $d\varphi$ is surjective and $d\varphi (ad(x)X)= ad(\varphi(x))d\varphi(X)$, we have $\nl(\varphi\gamma) = d\varphi\left(\nl(\gamma)\right)$ and $\mathfrak{p}(\varphi\gamma)=d\varphi \left(\mathfrak{p}(\gamma)\right)$. We conclude that $\varphi\gamma$ is the desired cocharacter for $X+\zl$.
\end{proof}

\begin{ste}
 Howe's conjecture holds in $PGL_n$.
\end{ste}
\begin{proof}
 We follow Harish-Chandra \cite{HC99} again and mention the adjustments. We replace $\cn$ by $N=\;^G\nl$. The proof of Harish-Chandra uses three properties of $\cn$ (in brackets the Lemma's in \cite{HC99} where the property is used):
 \begin{enumerate*}
  \item $\cn\cap S$ is compact (Lemma 11.9)
  \item For all compact subsets $\omega$ in $\gl$ there exists a lattice $L_1$, such that $^G\omega \subset L_1+\cn$. (Lemma 12.2)
  \item If $c\in \F$ and $Y\in \cn$, then $cY\in \cn$ (Lemma 12.3)
 \end{enumerate*}
 By Corollary \ref{NccN} (1) and (3) also hold for $N$ and (2) is Lemma \ref{L12.2}. Now we are left to prove Theorem 13.1 for $X\in N\cap S$. By Corollary \ref{gecax} we can use the proof of Theorem \ref{Higc}.
\end{proof}

\begin{lem}
 Let $G$ be a $\F$-split group. If $T^{ad}$ is a $\F$-split torus of $G^{ad}$, then $Ad^{-1}(T^{ad})$ is a $\F$-split torus of $G$.
\end{lem}
\begin{proof}
Without loss of generality we assume that $T^{ad}$ is a maximal $\F$-split torus of $G^{ad}$. Let $S$ be a maximal split torus of $G$ and $B$ a Borel group containing $S$. Then $S^{ad} := Ad(S)$ is a maximal split torus of $G^{ad}$ and $B^{ad}$ a Borel group containing $S^{ad}$. Take $g\in G^{ad}(\F)$ such that $gS^{ad}g^{-1} = T^{ad}$. Take $w^{ad}\in W^{ad}$ such that $g\in U_{(w^{ad})^{-1}}w^{ad}B^{ad}$. By multiplying $g$ with a suitable element of $S^{ad}$, we may assume that $g\in U_{(w^{ad})^{-1}}w^{ad}U^{ad}$. Take $w\in W$ such that $Ad(w)=w^{ad}$, then
\[ Ad : U_{w^{-1}}wU \rightarrow  U_{(w^{ad})^{-1}}w^{ad}U^{ad}\]
is a bijection. Therefore there exists a $h\in G(\F)$ such that $Ad(h)=g$. Thus \[Ad(hSh^{-1}) = gS^{ad}g^{-1}=T^{ad}.\]
Thus $Ad^{-1}(T^{ad})=hSh^{-1}$ is a $\F$-split torus.
\end{proof}

\begin{pro}
 Let $G$ be a $\F$-split group whose normal connected simple parts are all groups of type $A$. If $p{\not|}\kappa_v(G)$, then Howe's conjecture holds for $G$.
\end{pro}
\begin{proof}
Since $p{\not|}\kappa_v(G)$, the map $Ad : G \rightarrow G^{ad}$ is separable. Thus $d(Ad) : \gl\rightarrow \gl^{ad}$ is surjective.\\
Let $N := \{ n\in \gl \mid \exists(\lambda \in X_*(G)(\F)) \lim_{t\rightarrow 0} \lambda(t)n\lambda(t^{-1})=0\}$. Define $N^{ad}$ is the same way. Certainly $d(Ad)(N)\subset N^{ad}$.  We will show $N = \mathcal{N}\cap d(Ad)^{-1}(N^{ad})$. Let $n_a\in N^{ad}$. Take $\gamma_a \in X_*(G^{ad})(\F)$ such that
\[ \lim_{t\rightarrow 0} \gamma_a(t)n_a\gamma_a(t^{-1})=0.\]
Let $T^{ad}$ be a maximal $\F$-split torus that contains the image of $\gamma_a$.\\ Then $T := Ad^{-1}(T^{ad})$ is also a maximal $\F$-split torus. 

Let $\gamma\in X_*(T)$ and $m\in\N_{>0}$ be such that $Ad\circ \gamma = m\gamma_a$. Take a $n\in\mathcal{N}$ with $d(Ad)(n)=n_a$. Since $n\in \mathcal{N}$, $\gamma(t)n\gamma(t)^{-1} \in \mathcal{N}$ for all $t\in \F$. Because $\ker d(Ad) \subset \tl$, $d(Ad)$ restricted to $\mathcal{N}$ is a bijection. Now $ \lim_{t\rightarrow 0} \gamma_a^m(t)n_a\gamma_a^{-m}(t)=0$, thus
\[\lim_{t\rightarrow 0} \gamma(t)n\gamma(t)^{-1} =0.\] 
Therefore $n\in N$. Thus $N = \mathcal{N}\cap d(Ad)^{-1}(N^{ad})$.\\
The subset $N$ is closed, because $N_a$ is closed, $\mathcal{N}$ is closed and $d(Ad)$ is continuous. It is clear from the definition of $N$ that $cN=N$ for all $c\in \F^\times$.

Let $\gamma \in X_*(T)$ such that $Ad\circ\gamma= m\gamma$ for some $m\in \N_{>0}$.\\
Since $Ad$ is separable, $d(Ad) \gl(\geq 0,\gamma) = \gl^{ad}(\geq 0 ,\gamma_a)$. By the calculations on $PGL$ we know that
\[ [n_a,\mathfrak{p}^{ad}] = \nl^{ad}. \]
Since $d(Ad)$ is a bijection restricted to $\nl$,
\[ [n,\mathfrak{p}] = \nl.\]
Thus we can use the proof of Theorem \ref{Higc}.
\end{proof}

\subsection{The Howe's conjecture classification ($\F$-split case)}
In this subsection we determine exactly for which $\F$-split reductive groups Howe's conjecture holds. 
\begin{ste}\label{Hscc}
 Let $G$ be a reductive $\F$-split group, then the following statements are equivalent
 \begin{enumerate}
  \item The characteristic $p$ of $\F$ is good and $p{\not|}\kappa_v(G)$
  \item For all compact subsets $\omega$ and lattices $L$ in $\gl$:
  \[ \dim J_L(\omega) <\infty\]
 \end{enumerate}
\end{ste}
\begin{proof}
If the characteristic $p$ of $\F$ is bad, then $G$ has bad pairs. So in that case Howe's conjecture does not hold.\\
If $p|\kappa_v(G)$, then Howe's conjecture does not hold by Theorem \ref{scnhc}.\\
Assume that $p$ is good and $p{\not|}\kappa_v(G)$.\\
Let $G = R(G)G_1\cdots G_m$ with $G_i$ connected normal simple groups. Without loss of generality assume that $G_1,\ldots, G_k$ are the only groups of type $A$. Let $T$ be a maximal $\F$-split torus. Define \[G_A(T) := T \prod_{i=1}^k G_i.\]
Since $p{\not|}\#\text{coker}\;\Phi$, also $p{\not|}\#\text{coker} X_*(T)\rightarrow \text{Hom}_\Z(\Z R(G_A(T),T),\Z)$. Thus we can use Harish-Chandra's method for $G_A$. Because $p$ is good we can also use Harish-Chandra's method in the connected normal simple groups $G_i$ with $i>k$. Therefore we can use Harish-Chandra's method for the whole group: again we have to substitute $\mathcal{N}$ in the proof of Harish-Chandra.\\
\[ N := N_A\oplus \left(\bigoplus_{i={k+1}}^m \mathcal{N}_i\right),\]
where $N_A := \{ n\in\gl_A \mid \exists(\lambda \in X_*(G)(\F)) \lim_{t\rightarrow 0} \lambda(t)n\lambda(t^{-1})=0\}$.\\
We are left with proving the following about $N$:
 \begin{enumerate*}
  \item $N\cap S$ is compact (Lemma 11.9)
  \item For all compact subsets $\omega$ in $\gl$ there exists a lattice $L_1$, such that $^G\omega \subset L_1+N$. (Lemma 12.2)
  \item If $c\in \F$ and $Y\in N$, then $cY\in N$ (Lemma 12.3)
  \item For all $n\in N$ there exists a cocharacter $\gamma\in X_*(G)(\F)$ such that $\nl\subset [n,\gl]$. (\S 13.1)
 \end{enumerate*}
The statements for $N$ follow from the fact that they are true for $N_A$ and $\mathcal{N}_i$.
\end{proof}

\begin{gev}
 If $G$ is $\F$-split and has finitely many nilpotent orbits, then Howe's conjecture holds for $G$.
\end{gev}
\begin{proof}
If the characteristic $p$ of $\F$ is bad for $G$ or if $p|\kappa_v(G)$, then there are infinitely many nilpotent orbits.
\end{proof}

\section{The separable classification}
In this section we give a characterization of the reductive groups whose nilpotent orbits are all separable. As a consequence we get a large class of reductive groups for which the number of nilpotent orbits is finite and Howe's conjecture holds. We take a look at the cokernels of the following functions:
\begin{align*}
 \Phi &: X_*(T) \rightarrow \text{Hom}_\Z(\Z\Delta,\Z),\\
 \Phi^\vee &: X^*(T) \rightarrow \text{Hom}_\Z(\Z\Delta^\vee,\Z).
\end{align*}

\begin{lem}\label{lol}
$p|\rho_v(G)$ if and only if the $H_\alpha := d\alpha^\vee(1)\in\tl$, for $\alpha\in \Delta$, are linearly dependent.
\end{lem}
\begin{proof}
We have the following isomorphism of vector spaces: $\tl \cong X_*(T)\otimes_\Z \F$. Let $\eps_i$ be a basis for $X^*(T)$ and $\eps_i^\vee$ a dual basis in $X_*(T)$. Let $\alpha^\vee\in X_*(T)$. Now $\alpha^\vee = \sum_{i=1}^m \left<\eps_i,\alpha^\vee\right>\eps_i^\vee$. Hence $d\alpha^\vee(1) = \sum_{i=1}^m \left<\eps_i,\alpha^\vee\right>d\eps_i^\vee(1)$. Let $\alpha^\vee_1,\ldots,\alpha^\vee_n$ be the simple roots in $\Delta^\vee$. Define $M$ to be the $n\times m$ matrix with the following entries 
\[ M_{ij} := \left<\eps_j,\alpha^\vee_i\right>.\]
Then $M$ is the matrix corresponding with the map $\Phi^\vee$.\\
The matrix $M^{tr}$ is the matrix corresponding with the linear span of the $H_{\alpha_i}$'s.\\
Let $(d_1,\ldots,d_n)$ be the entries on the diagonal of the Smith normal form of $M$. Then $\rho_v(G)=\#\coker \Phi^\vee = \prod^n_{i=1}d_i$. The linear span of the $H_{\alpha_i}$'s is $n$-dimensional if and only if $p{\not|}\prod^n_{i=1}d_i$. 
\end{proof}

\begin{ste}\label{cnos}
The nilpotent orbits are separable if and only if the $p$ is good and $p {\not|} \kappa_v(G)$ and $p{\not|} \rho_v(G)$.
\end{ste}

\begin{proof}[Proof $\Rightarrow$.]
If $p$ is bad or divides $\kappa_v(G)$, then the regular nilpotent orbit is inseparable by Corollary \ref{rnons} and Theorem \ref{caisrn}. Assume that $p$ divides $\rho_v(G)$. Let $X := \sum_{\alpha\in\Delta}E_\alpha$. Then:
\[[\sum_{\alpha\in\Delta} E_\alpha,\sum_{\alpha\in-\Delta} c_\alpha E_\alpha]= \sum_{\alpha\in\Delta} c_\alpha H_\alpha.\] 
Now $p$ divides the cokernel exactly when the $H_\alpha = d\alpha^\vee(1)$ are linearly dependent. Thus there exists a $Y\in \nl_{-1}-\{0\}$, such that $[X,Y]=0$. Since $Z_G(X)\subset B$ and $\nl_{-1}\cap \bl=0$, the orbit of $X$ is not separable.
\end{proof}

Before we prove the implication in the other direction, we first state a few lemma's.

\begin{lem}\label{dsgl}
If $p$ is good for $G$, then
\[ \gl=\gl_A\oplus\gl_1\oplus\cdots\oplus \gl_m\]
with $G =G_A\prod_{i=1}^m G_i$ where $G_i$ are all the closed normal connected simple groups not of type $A$ and $G_A$ is generated by $R(G)$ and the closed normal connected simple groups of type $A$ in $G$.
\end{lem}
\begin{proof}
We have $Ad(G)\cong G^{ad}_A\times G_1^{ad}\times \cdots\times G_m^{ad}$.
Let $\Pi : Ad(G)\rightarrow G_1^{ad}\times \cdots \times G_m^{ad}$ be the corresponding projection map. Since $p$ is good for $G$, it is very good for $G_{c} = G_1\cdots G_m$. Thus the map $d(Ad) : \gl_{c}\rightarrow \gl^{ad}_{c}$ is surjective. Since $\dim G_{c}=\dim Ad(G_{c})$, it is a bijection. Thus $d(\Pi\circ Ad):\gl_{c}\rightarrow \gl_{c}^{ad}$ is a bijection. Therefore $\gl_{c} = \gl_1\oplus\cdots\oplus \gl_m$, since $p$ is very good for $G_{c}$. Moreover $\ker d(\Pi\circ Ad) \cap \gl_{c}=0$. Since $\gl_A\subset \ker d(\Pi\circ Ad)$, also $\gl_A\cap \gl_{c}=0$. Because $\dim \gl_A +\dim \gl_{c}=\dim \gl$, the Lemma follows.
\end{proof}

\begin{lem}
If $p$ is good for $G$ and $p{\not|}\kappa_v(G)$, then
$d(Ad) :\gl_A\rightarrow \gl_A^{ad}$ is surjective.
\end{lem}
\begin{proof}
By Lemma \ref{dsgl} and its proof we have $\gl = \gl_A\oplus \gl_c$ and $d(Ad):\gl_c\rightarrow \gl_c^{ad}$ is surjective. Since $p{\not|} \kappa_v(G)$ the map $d(Ad) : \gl \rightarrow \gl_A^{ad}\oplus \gl_c^{ad}$ is surjective. Let $\Pi_A : Ad(G)\rightarrow G_A^{ad}$, then $d(\Pi_A\circ Ad) : \gl \rightarrow \gl_A^{ad}$ is surjective. Since $\gl_c$ is contained in its kernel and $\gl = \gl_A\oplus \gl_c$, $d(Ad) : \gl_A\rightarrow \gl_A^{ad}$ is surjective.
\end{proof}
\begin{gev}\label{pndca}
 If $p$ is good for $G$ and $p{\not|}\kappa_v(G)$, then $p$ does not divide the cokernel of the following map:
 \[ \Phi_A : X_*(T_A)\rightarrow \hom_Z(\Z\Delta_A,\Z).\]
\end{gev}
\begin{proof}
 The map $d(Ad) : \gl_A\rightarrow \gl_A^{ad}$ is surjective, thus  $p{\not|}\#\coker \Phi_A$ by Proposition \ref{csad}. 
\end{proof}

\begin{lem}
 Let $n = n_1+\cdots+n_m$ and $\gamma_i$ be the cocharacters $\gamma_i\in X_*(T\cap G_i^{ad})$ associated with $n_i$ in $\gl_i^{ad}$. Let $\gamma\in X_*(T)$ be the cocharacter associated with $n$ in $G$. Then $d(Ad)\circ \gamma= \sum_{i=1}^m \gamma_i$
\end{lem}
\begin{proof}
 Clearly $\sum_{i=1}^m \gamma_i$ is a cocharacter associated with $n$ in $\gl^{ad}$. Also $d(Ad)\circ \gamma$ is a cocharacter associated with $n$ in $\gl^{ad}$. (See \cite[\S5.6]{Ja04}) Since there is at most one cocharacter of $Ad(T)$ associated with $n$ by \cite[Corollary 22]{Mc04}, they are equal.
\end{proof}

\begin{lem}
 Let $G = GL_m$ and $G^{ad} =PGL_m$. Let $n\in \gl^{ad}$ be a nilpotent element with associated cocharacter $\gamma$. Then
 \begin{align*}
 &[n,\gl^{ad}(k)]=\gl^{ad}(k+2) &\text{ for } &k\geq -1\\
 &[n,\cdot]:\gl^{ad}(k)\rightarrow \gl^{ad}(k+2) \text{ is injective}&\text{ for }&k=-1\text{ and }k\leq -3
 \end{align*}
\end{lem}
\begin{proof}
For $GL_m$ and $n\in\gl_m$ nilpotent $[n,\gl(k)] = \gl(k+2)$ for $k\geq -1$ and $[n,\cdot] : \gl(k)\rightarrow \gl(k+2)$ is injective for $k\leq -1$. Since the map \[d(Ad):\bigoplus_{k\geq 1}\;\gl(-k)\oplus \gl(k)\rightarrow \bigoplus_{k\geq 1}\;\gl^{ad}(-k)\oplus \gl^{ad}(k)\]
is a bijection, $[n,\gl^{ad}(k)]=\gl^{ad}(k+2)$ for $k\geq -1$ and $[n,\cdot]:\gl^{ad}(k)\rightarrow \gl^{ad}(k+2)$ is injective for $k=1$ and $k\geq -3$.
\end{proof}

\begin{lem}\label{pcne}
Let $G=GL_m$. Let $n = \sum_{\alpha\in \Gamma} c_\alpha n_\alpha$, with $\Gamma\subset \Delta$. If $[n,m]\in\zl$ and $m\in \gl(-2)$, then $m=\sum_{\alpha\in-\Gamma} d_\alpha n_\alpha$ for some $d_\alpha\in \F$. 
\end{lem}
\begin{proof}
Let $\gamma\in X_*(T)$ be such that there exists a $l\in \N_{\geq 0}$ such that for all $\alpha\in \Delta$: 
\[\left<\alpha,\gamma\right>=\left\{ \begin{array}{cl} l &\text{if } \alpha\in \Gamma,\\ 0 &\text{if } \alpha\not\in \Gamma.\end{array}\right.\]
We know that $[n,\cdot]: \gl(-2;\tau)\rightarrow \gl(0;\tau)$ is injective for the associated cocharacter $\tau\in X_*(T)$ of $n$. Define $\gl(-2) := \gl(-2;\tau)$ and $\gl_i(-2) := \gl(-2;\tau)\cap \gl(il;\gamma)$. Then
\[ \gl(-2) = \bigoplus \gl_i(-2)\]
and
\[ [n,\gl_i(-2)] \subset \gl(0;\tau)\cap \gl(l(i+1);\gamma).\]
Because $\zl \subset \gl(0;\gamma)$ and $[n,\cdot]|_{\gl(-2)}$ is injective, then $m\in \gl(-l;\gamma)$. 
\end{proof}

\begin{proof}[Proof Theorem \ref{cnos} $\Leftarrow$.]
Let $n\in \nl$. Take $n_A\in\nl_A$ and $n_i\in\nl_i$, such that
\[ n = n_A+n_1+\cdots+n_m.\]
Then
\[ Z_G(n) = Z_{G_A}(n_A)\prod_{i=1}^m Z_{G_i}(n_i) \]
\[ Z_\gl(n) = Z_{\gl_A}(n_A)\oplus\bigoplus_{i=1}^m Z_{\gl_i}(n_i)\]

Since the $G_i$ are simple and $p$ is very good for $G_i$, the $G_i$-orbit of $n_i$ is separable:
\[\dim Z_{\gl_i}(n_i)=\dim Z_{G_i}(n_i).\]

Thus we are left with showing that $\dim Z_{G_A}(n_A) = \dim Z_{\gl_A}(n_A)$. Since $p$ is good for $G$, it is also good for $G_A$. By Corollary \ref{pndca} and Lemma \ref{lol}, $p$ does not divide the order of the cokernels of $\Phi_A$ and $\Phi_A^\vee$.

Thus without loss of generality we assume that $G$ only consists of groups of type $A$ and a center. Thus $G^{ad} = \prod_{i=1}^k PGL_{n_i}$. 

Since $p\;{\not|}\kappa_v(G)$ the map $d(Ad):\gl\rightarrow \gl^{ad}$ is surjective. Let $n\in\gl$ be nilpotent and $\gamma$ be a cocharacter associated with $n$. Define $P := P(\gamma)$. Then $Ad\circ \gamma$ is a cocharacter associated with $Ad(n)$. For $G^{ad}$ the following holds:
\begin{align*}
 &[n,\gl^{ad}(k)]=\gl^{ad}(k+2) &\text{ for } &k\geq -1\\
 &[n,\cdot]:\gl^{ad}(k)\rightarrow \gl^{ad}(k+2) \text{ is injective}&\text{ for }&k=-1\text{ and }k\leq -3
\end{align*}
Since $d(Ad)$ is surjective and injective on the nilpotent elements, then
\begin{align*}
 [n,\pl]&=\gl(\geq 2)\\
 [n,\cdot]&:\gl(k)\rightarrow \gl(k+2) \text{ is injective for }k=-1\text{ and }k\leq -3.
\end{align*}
Therefore 
\[\dim Z_G(n) = \dim Z_P(n) = \dim Z_{\pl}(n),\]
since $\overline{Ad\;P(n)}=\gl(\geq 2)$ and $[n,\pl]=\gl(\geq 2)$.\\
If $Z_\gl(n)\cap \gl(k)=0$ for $k\leq -1$, then $Z_\gl(n)=Z_\pl(n)$.\\
For $k=-1$ and $k\leq -3$ the function $[n,\cdot]:\gl(k)\rightarrow \gl(k+2)$ is injective. Thus $Z_\gl(n)\cap \gl(k)=0$, for $k=-1$ and $k\leq -3$.\\
Thus we are done if the kernel of $[n,\cdot] :\gl(-2)\rightarrow \gl(0)$ is $0$.

In $G^{ad}$ every nilpotent element is conjugated with an element of the form $\sum_{\alpha\in \Gamma}E_\alpha$, with $\Gamma\subset \Delta$.
Let $n = \sum_{\alpha\in \Gamma} E_\alpha$ with $\Gamma\subset \Delta$ and $m\in \gl(-2)$. If $[n,m]=0$, then $[d(Ad)(n),d(Ad)(m)]=0$. Thus $m = \sum_{\alpha\in-\Gamma} c_\alpha E_\alpha$ for some $c_\alpha\in \F$ by Lemma \ref{pcne}. Now
\[0=[n,m] = \sum_{\alpha\in \Delta} c_\alpha H_\alpha.\]
Because $p\;{\not|}\rho_v(G)$, the $H_\alpha$ are linearly independent. Thus $c_\alpha = 0$ for all $\alpha \in -\Gamma$, hence $m=0$.\\
Thus every nilpotent orbit is separable.
\end{proof}

\section{On the number of nilpotent orbits}
In this section we discuss when the number of nilpotent orbits is finite.

\begin{ste}\cite[Theorem 40]{Mc04}\label{s40mc}
If $p$ is good and all the nilpotent orbits are separable, then there are only finitely many nilpotent orbits.
\end{ste}
\begin{gev}\label{ecmc}
If $p$ is good and $p {\not|} \kappa_v(G)$ and $p{\not|} \rho_v(G)$, then there are only finitely many nilpotent orbits.
\end{gev}
\begin{proof}
The condition in the Corollary is equivalent to the one in Theorem \ref{s40mc} by Theorem \ref{cnos}.
\end{proof}
In this section we will prove the converse of Corollary \ref{ecmc}. If $G$ is $\F$-split and $p$ is bad or divides $\kappa_v(G)$, then there are infinitely many regular nilpotent orbits by Theorem \ref{bnhc} and Proposition \ref{ubaan}. So it is enough to prove that if $G$ is $\F$-split, $p$ is good, $p{\not|}\kappa_v(G)$ and $p|\rho_v(G)$, then $G$ has infinitely many nilpotent orbits. First a Theorem that we can easily deduce from the theory of the previous section.

\begin{ste}\label{fnivg}
 If $G$ is semi-simple and the characteristic of $\F$ is not very good, then there are infinitely many nilpotent orbits.
\end{ste}
\begin{proof}
If the characteristic of $\F$ is bad, then we have already showed that there are infinitely many nilpotent orbits. So without loss of generality we assume  $G$ is has at least one normal simple groups of type $A_n$, with $p|n+1$. Now the proof is split in two cases: $p\; | \kappa_v(G)$ and $p\; {\not|} \kappa_v(G)$.\\

If $p\; | \kappa_v(G)$, then $\gl$ has infinitely many nilpotent orbits by Proposition \ref{ubaan}.\\

If $p\;{\not|}\kappa_v(G)$, then $d(Ad) : \gl\rightarrow \gl^{ad}$ is an isomorphism by Theorem \ref{csad}. Since there are infinitely many nilpotent orbits in $\gl^{ad}$ by Corollary \ref{pglio}, there are also infinitely many nilpotent orbits in $\gl$. 
\end{proof}

\begin{pro}\label{ppddc}
 Let $G$ be a reductive group with only normal simple subgroups of type $A$ for which $p$ is not very good. Assume that $p|\rho_v(G)$. Let $H$ be a reductive group with $G\nd H$. Let $\mathcal{N}$ be the set of nilpotent elements of $\gl$. Then there are infinitely many nilpotent $H$-orbits in $\mathcal{N}$. 
\end{pro}
\begin{proof}
The proof of this proposition is distributed over two lemmas.
\begin{lem}
 If $\alpha^\vee\in \Delta^\vee$, then $Ad\circ \alpha^\vee \in \Delta_{ad}^\vee$.
\end{lem}
\begin{proof}
 The reader could verify this by taking the Chevalley basis on $\gl$.
\end{proof}

Let $\Delta^\vee = \{ \alpha^\vee_{11},\ldots \alpha^\vee_{nm_{n}}\}$, such that $\alpha^\vee_{ij}$ is connected in the Dynkin diagram with $\alpha^\vee_{i'j'}$ if and only if $i=i'$ and $j=j+1$.

\begin{lem}\label{tcl}
If $\sum_{\alpha\in\Delta}c_\alpha d\alpha^\vee(1) =0$, then for every $i$ there exists a $c_i$ such that $c_{\alpha_{ij}} = jc_i$ for all $j$. 
\end{lem}
\begin{proof}
Because $\sum_{\alpha\in\Delta}c_\alpha d\alpha^\vee(1) =0$, also
\[\sum_{\alpha\in\Delta_{ad}}c_\alpha d(Ad\circ\alpha^\vee)(1) =0.\]
( $d(Ad)\left(d\alpha^\vee(1)\right)=d(Ad\circ \alpha^\vee)(1)$ )\\
Since $\gl^{ad} = \gl_1\oplus \cdots \oplus \gl_n$, with $\gl_i$ the Lie algebra of $PGL_{m_i+1}$, then for every $i$:
\[\sum_{j=1}^{m_i} c_{\alpha_{ij}} d(Ad\circ \alpha_{ij}^\vee)(1)=0.\]
A small calculation in $\gl_i$ shows that there exists a $c_i$ such that $c_{\alpha_{ij}}=jc_i$.
\end{proof}

Since $p|\rho_v(G)$, by Lemma \ref{tcl} there exist $c_i\in \F$ such that
\[\sum_{i=1}^n\sum_{j=1}^{m_i} c_ijd\alpha_{ij}^\vee(1)=0\]
and at least one of the $c_i\not=0$. Without loss of generality assume that $1,\ldots,k$ are the $i$ with $c_i\not=0$.\\

Let $i\leq k$. Let $M_i(x)$ be the block matrix consisting of $\frac{m_i}{p}$ blocks of $p\times p$-matrices, with on each block the following matrix
\[ \left(\begin{array}{cccc} 0&\cdots&0&c_ix\\1&0&\cdots&0\\0&\ddots&\ddots&\vdots\\ \ddots&0&1&0\end{array}\right).\]
Thus the entries of $M_i$ are as follows:
\[ (M_i)_{kl} := \left\{\begin{array}{cl} 1 &\text{if } k=l+1 \text{ and } p{\not|}l\\
c_ix &\text{if } l=k+p-1\text{ and } p|l\\
0&\text{otherwise}
\end{array}\right.\]

Then $M_i(x)^p = c_ixI_{n_i}$.\\
Let $N(x)$ be the element in $\gl$ corresponding with $M_1(x)\oplus \cdots \oplus M_k(x)$. Then 
\[N(x)^p = \sum_{i=1}^n\sum_{j=1}^{m_i} xc_ijd\alpha_{ij}^\vee(1)=0.\]
Thus $N(x)$ is nilpotent.\\
Let $q$ be a power of $p$ such that $\mathcal{N}^{ad} = \{ X\in\gl^{ad} \mid X^q=0\}$.\\
Let $\phi' : \mathcal{N}^{ad} \rightarrow \F/\F^q$ be the following function:\\
Take $X_i\in \Gl_{m_i+1}$ such that $X = \oplus_{i=1}^{n} X_i+\zl_i$. Then for each $i$ we have a $z_i$ such that $X_i^q=z_iI_{m_i+1}$ in $\Gl_{m_i+1}$. Define $\phi'(X) := z_1$.\\
If $X'_i$ are also representatives for $X$, then $z'_i=z_i+a_i^q$ for $a_i\in \F$. Thus $\phi$ is well-defined.\\
Since $G^{ad} = \prod_{i=1}^n PGL_{m_i+1}$, $\phi$ is also $G^{ad}$-invariant. Define $\phi : \mathcal{N}\rightarrow \F/\F^q$ by $\phi := \phi'\circ d(Ad)$. Then $\phi$ is $H$-invariant, since $H$ acts on $\mathcal{N}$ by conjugation and $G_1\nd H$. For $x\in \F$,
\[ \phi(N(x)) = x^{\frac{q}{p}}.\]

Since $\F^{\frac{q}{p}}/\F^q \cong \F/\F^{(p)}$ is infinite and $\phi$ is $H$-invariant, there are infinitely many nilpotent $H$-orbits.
\end{proof}

\begin{ste}\label{tpddc}
 If $p$ is good and $p{\not|}\kappa_v(G)$, but $p|\rho_v(G)$, then there are infinitely many nilpotent orbits.
\end{ste}
\begin{proof}
Let $G=R(G)G_1\cdots G_l$ with $G_i$ the minimal simple normal connected subgroups of $G$. Assume that $G_1,\ldots,G_n$ are the groups of type $A$ for which $p$ is not very good. Define $G_A := R(G)G_1\cdots G_n$ and $G_C := G_{n+1}\cdots G_l$. Because $p{\not|}\kappa_v(G)$, the map $d(Ad) : \gl\rightarrow \gl^{ad}$ is surjective. Because $p$ is very good for $G_C$, $d(Ad) : \gl_C\rightarrow \gl_C^{ad}$ is surjective. Since $G_C$ is semi-simple, the map is even an isomorphism.
\[\gl^{ad} \cong \bigoplus_{i=1}^l \gl^{ad}_i = \gl^{ad}_A\oplus \gl^{ad}_C\]
Define $Ad_C: \gl\rightarrow \gl^{ad}_C$ by the composition of the projection and $d(Ad) :\gl\rightarrow \gl^{ad}$. Then $\gl_A\subset \ker Ad_C$ and $\ker Ad_C \cap \gl_C =0$. Hence $\gl_A\cap \gl_C=0$, thus $\gl = \gl_A\oplus \gl_C$. By Lemma \ref{lol} the $d\alpha^\vee(1)$'s are linearly dependent. Because of the decomposition of $\gl$, the $d\alpha^\vee(1) : \alpha\in \Delta_A$ are linearly dependent or the $d\alpha^\vee(1) : \alpha\in \Delta_C$ are linearly dependent. Since $p$ is very good for $G_C$ the $d\alpha^\vee(1) : \alpha\in \Delta_C$ are linearly independent. So the $d\alpha^\vee(1): \alpha\in \Delta_A$ are linearly dependent. Therefore we can apply Proposition \ref{ppddc} with $H=G$ and $G=G_A$.
\end{proof}

\begin{ste}\label{cfnisc}
 If $G$ is $\F$-split, then the following are equivalent:
 \begin{enumerate*}
  \item The number of nilpotent orbits is finite
  \item All the nilpotent orbits are separable
  \item The regular nilpotent orbit is separable
  \item $p$ is good and $p{\not|}\kappa_v(G)\rho_v(G)$
 \end{enumerate*}
\end{ste}
\begin{proof}
 (2) implies (1) by \cite[Theorem 40]{Mc04}. (1) implies (4) by Theorem \ref{tpddc}, Theorem \ref{bnhc} and Proposition \ref{ubaan}. (4) implies (2) by Theorem \ref{cnos}. By the proof of Theorem \ref{cnos} also not (4) implies not (3).
\end{proof}

\newpage
\section{Appendix}
\addtocontents{toc}{\protect\setcounter{tocdepth}{-1}}
In this appendix we will prove the following theorem for all the adjoint simple groups.
\begin{ste}
 Let $G$ be an adjoint simple group. Let $X := \sum_{\alpha\in\Delta}E_\alpha$. Assume $p= \text{char}\;\F$ is bad for $G$. Then
 \begin{enumerate*}
  \item $[X,\cdot] : \nl_i\rightarrow \nl_{i+1}$ is surjective for $1\leq i\leq p-1$
  \item $\dim \nl_{p+1}/[X,\nl_p] = 1$
  \item $|\Delta| = \dim \nl_1 = \dim \nl_{i} +1$ for $2\leq i\leq p+1$
 \end{enumerate*}
\end{ste}
The theorem will be proved with a case by case consideration of the adjoint simple groups. For the simply laced root systems we introduce the following notation. For all $i\in \Delta$ we choice a non-zero $E_i\in\gl_i$. For $i,j,k,l\in \Delta$ we define
\begin{align*}
 E_{kl} &:= [E_k,E_l]\\
 E_{jkl} &:= [E_j,E_{kl}]\\
 E_{ijkl} &:= [E_i,E_{jkl}]
\end{align*}
If $k+l\in R$ ($j+k+l\in R$, $i+j+k+l\in R$), then $E_{kl}\not=0$ ($E_{jkl}\not=0$, $E_{ijkl}\not=0$ respectively).

\subsection{$D_n$, $n\geq  4$}
For this part char $\F=2$.\\
We start with some calculations on $D_4$.
\[
 \xymatrix{ & 1 \ar@{-}[d]& \\ 2\ar@{-}[r]&4\ar@{-}[r]&3}
\]
Thus $\Delta = \{1,2,3,4\}$.
\begin{align*}
 \nl_1 &=\left<E_1,E_2,E_3,E_4\right>\\
 \nl_2 &= \left<E_{14},E_{24},E_{34}\right>\\
 \nl_3 &= \left<E_{124},E_{134},E_{234}\right>
\end{align*}
By looking in $GL_4$ one can show that it is possible to choice the basis as follows:
\begin{align*}
 [E_{i},E_{4}] &= E_{i4}\\
 [E_1,E_{i4}] &= E_{1i4}\\
 [E_2,E_{34}] &= E_{234}
\end{align*}
Then
\begin{align*}
 [E_2,E_{14}] &= [E_2,[E_1,E_4]]= -[E_1,[E_4,E_2]]-[E_4,[E_2,E_1]] = [E_1,[E_2,E_4] = E_{124}\\
 [E_3,E_{14}] &= [E_3,[E_1,E_4]]= -[E_1,[E_4,E_3]]-[E_4,[E_3,E_1]] = [E_1,[E_3,E_4] = E_{134}\\
 [E_3,E_{24}] &= [E_3,[E_2,E_4]]= -[E_2,[E_4,E_3]]-[E_4,[E_3,E_2]] = [E_2,[E_3,E_4] = E_{234}
\end{align*}
Let $X=E_1+E_2+E_3+E_4$.\\
Thus with these basses the action $[X,\cdot]:\nl_2\rightarrow \nl_3$ has the following matrix:
\[\left(\begin{array}{ccc} 1&1&0\\1&0&1\\0&1&1 \end{array}\right).\]
One sees that this matrix has at least rank $2$. Because the determinant is $-2=0$, the matrix has rank $2$.\\

Now for general $n\geq 4$. We take a slightly different numbering of the roots.
\[
 \xymatrix{ & n \ar@{-}[d]& & & &\\ 1\ar@{-}[r]&2\ar@{-}[r]&3 \ar@{--}[r] & &\ar@{--}[r] & n-1 }
\]
\begin{align*}
 \nl_2 &= \left< E_{2n},E_{i,i+1} : 1\leq i\leq n-2\right>\\
 \nl_2' &= \left< E_{i,i+1} : 1\leq i\leq n-2\right>\\
 \nl_3 &= \left<E_{12n},E_{23n}, E_{i,i+1,i+2} : 1\leq i\leq n-3\right>\\
 \nl'_3 &= \left<E_{i,i+1,i+2} : 1\leq i\leq n-3\right>
\end{align*}
Let $\pi : \nl_3\rightarrow \nl_3'$ be the projection on the basis. By looking at $GL_{n}$ (the roots $1,\cdots,n-1$ correspond with a group of type $A_n$) one sees that $\pi\circ \nl'_2\rightarrow \nl_3'$ is surjective. Also $[X,E_{2n}]\not=0$ lies in the vector space $\left<E_{12n},E_{23n}\right>$. Thus the rank of the $[X,\cdot]$ is at least $n-1$. By the calculations on $D_4$ we know that the kernel is not trivial. Therefore $[X,\cdot] :\nl_2\rightarrow \nl_3$ has a $1$-dimensional cokernel. 

\subsection{$B_n$}
For this part char $\F=2$.\\
We have seen the Theorem for $B_2$ in $SO_5$.\\
Now we look at $B_3$, the $B_n$ with $n\geq 4$ follow in the same way as $D_n$ is a consequence of $D_4$.\\
Let $\gl$ be the Lie algebra of $D_4$ and $\sigma$ be the action on $\gl$ corresponding with the permutation $(12)$ of the Dynkin diagram of $D_4$.
\begin{align*}
 \nl_1^\sigma &= \left<E_1+E_2,E_3,E_4\right>\\
 \nl_2^\sigma &= \left<E_{14}+E_{24},E_{34}\right>\\
 \nl_3^\sigma &= \left<E_{124},E_{134}+E_{234}\right>
\end{align*}
From the calculations on $D_4$ we know that the corresponding matrix is
\begin{align*}
 \left(\begin{array}{ccc} 1&1&0\\1&0&1\\0&1&1\end{array}\right).&\\
 \intertext{Thus the matrix for $[X,\cdot]:\nl^\sigma_2\rightarrow \nl_3$ is:}
 \left(\begin{array}{ccc} 1&1&0\\1&0&1\\0&1&1\end{array}\right)\left(\begin{array}{cc} 1&0\\1&0\\0&1\end{array}\right)&=\left(\begin{array}{cc}2&0\\1&1\\1&1\end{array}\right).\\
 \intertext{Therefore the matrix for $[X,\cdot]:\nl^\sigma_2\rightarrow \nl^\sigma_3$ is:}
 \left(\begin{array}{cc} 2&0\\1&1\end{array}\right).&
\end{align*}
This matrix has determinant $2=0$ and at least rank $1$. Thus the cokernel of $[X,\cdot]: \nl_2\rightarrow \nl_3$ is $1$ dimensional.
\subsection{$G_2$}
Let $\gl$ be the Lie algebra of $D_4$ and $\sigma$ by the action on $\gl$ corresponding with the permutation $(123)$ of the Dynkin diagram of $D_4$.\\
Let $E_{1234}$ be the basis for $\nl_4$. Now
\begin{align*}
 \nl_1^\sigma &= \left<E_1+E_2+E_3,E_4\right>\\
 \nl_2^\sigma &= \left<E_{14}+E_{24}+E_{34}\right>\\
 \nl_3^\sigma &= \left<E_{124}+E_{134}+E_{234}\right>\\
 \nl_4^\sigma &= \left<E_{1234}\right>
\end{align*}
The map $[X,\cdot]:\nl_2^\sigma\rightarrow \nl_3^\sigma$ is $(2)$ by a similar calculation as the $B_3$ case. Thus if char $\F=2$, then  $[X,\cdot]:\nl_2^\sigma\rightarrow \nl_3^\sigma$ has a $1$ dimensional cokernel.\\
The map $[X,\cdot]: \nl_3\rightarrow \nl_4$ is on the bases $\left(\begin{array}{ccc} 1&1&1 \end{array}\right)$. Therefore the map $[X,\cdot]:\nl_3^\sigma\rightarrow \nl_4^\sigma$ is $(3)$ with respect to the choice of the bases. Thus if char $\F=3$, then  $[X,\cdot]:\nl_3^\sigma\rightarrow \nl_4^\sigma$ has a $1$ dimensional cokernel.

\subsection{$C_n$}
In this subsection char $\F=2$.\\
For $C_2$, see $B_2$. Again we only have to look at $C_3$, because the Theorem for $C_n$, $n\geq 4$, follow in a similar fashion as for the groups of type $D_n$.\\
Number the simple roots of $A_5$ as follows:
\[ \xymatrix{ 1\ar@{-}[r]&2\ar@{-}[r]&3\ar@{-}[r]&4\ar@{-}[r]&5  } \]
Let $\gl$ be the adjoint Lie algebra of type $A_5$ and $\sigma$ be the action on $\gl$ corresponding with the permutation $(15)(24)$ of the Dynkin diagram of $A_5$. Then $\gl^\sigma$ is the adjoint Lie algebra of type $C_3$. We choose the following bases:
\begin{align*}
 \nl_1 &= \left<E_1,E_2,E_3,E_4,E_5\right> &\nl_1^\sigma&= \left<E_1+E_5,E_2+E_4,E_3\right>\\
 \nl_2 &= \left<E_{12},E_{24},E_{34},E_{45}\right> &\nl_2^\sigma&= \left<E_{12}-E_{45},-E_{23}+E_{34}\right>\\
 \nl_3 &= \left<E_{123},E_{234},E_{345}\right> &\nl_3^\sigma&=\left<E_{123}+E_{345},E_{234}\right>
\end{align*}
The matrix corresponding with $[X,\cdot] : \nl_2\rightarrow \nl_3$ is
\[ \left(\begin{array}{cccc} -1&1&0&0\\0&-1&1&0\\0&0&-1&1\end{array}\right).\]
Therefore the matrix for $[X,\cdot] : \nl^\sigma_2\rightarrow \nl^\sigma_3$ is
\[\left(\begin{array}{cc} -1&-1\\0&2\end{array}\right).\]
This matrix has determinant $-2=0$ and at least rank $1$. Thus the dimension of the cokernel of $[X,\cdot] : \nl^\sigma_2 \rightarrow \nl^\sigma_3$ is $1$.

\subsection{$E_6,E_7,E_8$}
In this part $G$ is a group of type $E_n$. We number the roots (if they exists) as follows
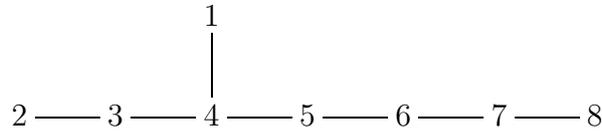
\[
 \xymatrix{ & & 1\ar@{-}[d] \\ 2\ar@{-}[r]&3\ar@{-}[r]&4\ar@{-}[r]&5\ar@{-}[r]&6\ar@{-}[r]&7\ar@{-}[r]&8 }
\]
$X := \sum_{i=1}^nE_i$
\subsubsection{char $\F=2$}
\begin{align*}
 \nl_2 &=\left<E_{14},E_{i,i+1} : 2\leq i\leq n-1\right>\\
 \nl'_2 &= \left<E_{i,i+1} : 2\leq i\leq n-1\right>\\
 \nl_3 &= \left<E_{134},E_{145}, E_{i,i+1,i+2} : 2\leq i \leq n-2\right>\\
 \nl'_3 &= \left< E_{i,i+1,i+2} : 2\leq i \leq n-2\right>
\end{align*}
Let $\pi : \nl_3\rightarrow \nl'_3$ be the projection. By looking at $GL_n$ we see that $[X,\cdot] : \nl'_2\rightarrow \nl'_3$ is surjective. Also $[X,E_{2n}]\not=0$ lies in the vector space $\left<E_{12n},E_{23n}\right>$. Thus the rank of the $[X,\cdot]$ is at least $n-1$. By the calculations on $D_4$ we know that the kernel is not trivial. Therefore $[X,\cdot] :\nl_2\rightarrow \nl_3$ has a $1$-dimensional cokernel.

\subsubsection{char $\F=3$}
First $E_6$, the others follow in a simular way as $D_n$ follows from $D_4$ in the characteristic $2$ case.
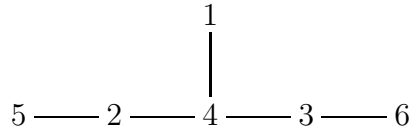
\[
 \xymatrix{ & & 1\ar@{-}[d] \\ 5\ar@{-}[r]&2\ar@{-}[r]&4\ar@{-}[r]&3\ar@{-}[r]&6 }
\]

\begin{align*}
\nl_2 &= \left<E_{14},E_{24},E_{34},E_{25},E_{36}\right>\\
\nl_3 &= \left<E_{124},E_{134},E_{234},E_{436},E_{425}\right>\\
\nl_4 &= \left<E_{1234},E_{1436},E_{1425},E_{3425},E_{2436}\right>
\end{align*}

The corresponding matrix of $[X,\cdot] : \nl_3\rightarrow \nl_4$ is
\[ \left(\begin{array}{ccccc} 1&1&1&0&0\\ 0&1&0&1&0\\ 1&0&0&0&1\\ 0&0&1&0&1\\ 0&0&1&1&0\end{array}\right).\]
The determinant of the matrix is $3=0$. Since the matrix has at least rank $4$, the dimension of cokernel of $[X,\cdot]:\nl_3\rightarrow \nl_4$ is $1$.

\subsubsection{char $\F=5$}
We use the construction of the adjoint Lie algebra of type $E_8$ as described in \cite[\S10.2]{Sp98} and adopt its notation. We define the bi-additive function $f$ on $X^*(T)$ as follows:
\[
 f(i,j) := \left\{\begin{array}{rl} -1 & \text{if } i \text{ is connected with } j \text{ and } i<j\\
                   1 & \text{if } i=j\\
                   0 & \text{otherwise}
                  \end{array}\right.
\]
If $\alpha,\beta \in R^+$, then $c_{\alpha,\beta} = (-1)^{f(\alpha,\beta)}$ if $\alpha+\beta\in R^+$ and $c_{\alpha,\beta}=0$ otherwise. Now $[e_\alpha,e_\beta]=c_{\alpha,\beta}e_{\alpha+\beta}$. Using this we can now construct the matrix corresponding with the function $[X,\cdot]$ with respect to the following bases: 
\begin{align*}
 \nl_5 &=\left<e_{12345},e_{13456},e_{14567},e_{23456}, e_{34567},e_{45678},e_{13445}\right>,\\
 \nl_6 &= \left<e_{123456},e_{123445},e_{134456},e_{134567},e_{145678},e_{234567},e_{345678} \right>,
\end{align*}
where $e_{a_1\cdots a_k}$ is short for $e_{b}$ with $b$ the root equal to the sum of the roots $a_i$.\\
The corresponding matrix is
\[
 \left(\begin{array}{ccccccc}
 1&-1&0&-1&0&0&0\\
 1&0&0&0&0&0&-1\\
 0&1&0&0&0&0&1\\
 0&1&-1&0&-1&0&0\\
 0&0&1&0&0&-1&0\\
 0&0&0&1&-1&0&0\\
 0&0&0&0&1&-1&0
\end{array} \right).
\]
It has determinant $5=0$ and rank $6$. Thus the cokernel of $[X,\cdot]:\nl_5\rightarrow \nl_6$ is $1$ dimensional.

\subsection{$F_4$}
The group $F_4$ is a folding of $E_6$. Let $\sigma$ be the action on the Lie algebra $\gl$ of $E_6$ corresponding with the isomorphism $(23)(56)$ of the Dynkin diagram of $E_6$. The Lie algebra of $F_4$ is $\gl^\sigma$.
\begin{align*}
 \nl^\sigma_1 &= \left<E_1,E_4,E_2+E_3,E_5+E_6\right>\\
 \nl_2^\sigma &= \left<E_{14},E_{24}+E_{34},E_{25}+E_{36}\right>\\
 \nl^\sigma_3 &= \left<E_{124}+E_{134},E_{234},E_{436}+E_{425}\right>\\
 \nl^\sigma_4 &= \left<E_{1234},E_{1436}+E_{1452},E_{3425}+E_{2436} \right>
\end{align*}

The matrix corresponding with $[X,\cdot]:\nl_2\rightarrow \nl_3$ is
\begin{align*} \left(\begin{array}{ccccc} 1&1&0&0&0\\ 1&0&1&0&0\\ 0&1&1&0&0\\ 0&0&1&0&1\\ 0&1&0&1&0\end{array}\right).&\\
\intertext{Thus the matrix of $[X,\cdot]:\nl_2^\sigma\rightarrow \nl_3$ is equal to:}
\left(\begin{array}{ccccc} 1&1&0&0&0\\ 1&0&1&0&0\\ 0&1&1&0&0\\ 0&0&1&0&1\\ 0&1&0&1&0\end{array}\right)\left(\begin{array}{ccccc}1&0&0&0&0 \\ 0&1&0&0&0 \\ 0&1&0&0&0 \\ 0&0&1&0&0 \\ 0&0&1&0&0  \end{array}\right)&=\left(\begin{array}{ccccc} 1&1&0&0&0\\1&1&0&0&0\\0&2&0&0&0\\0&1&1&0&0\\0&1&1&0&0\end{array}\right).
\intertext{Therefore the matrix of $[X,\cdot]: \nl_2^\sigma\rightarrow \nl_3^\sigma$ is}
\left(\begin{array}{ccc}1&1&0\\0&2&0\\0&1&1 \end{array}\right).&
\end{align*}
The determinant of this matrix is $2$. Thus if char $\F=2$, then the dimension of the cokernel of the map $[X,\cdot] :\nl_2^\sigma\rightarrow \nl_3^\sigma$ is $1$.\\

The matrix corresponding with $[X,\cdot]:\nl_3\rightarrow \nl_4$ is
\begin{align*} \left(\begin{array}{ccccc} 1&1&1&0&0\\ 0&1&0&1&0\\ 1&0&0&0&1\\ 0&0&1&0&1\\ 0&0&1&1&0\end{array}\right).&\\
\intertext{Thus the matrix of $[X,\cdot]:\nl_3^\sigma\rightarrow \nl_4$ is equal to:}
\left(\begin{array}{ccccc} 1&1&1&0&0\\ 0&1&0&1&0\\ 1&0&0&0&1\\ 0&0&1&0&1\\ 0&0&1&1&0\end{array}\right)\left(\begin{array}{ccccc}1&0&0&0&0 \\ 1&0&0&0&0 \\ 0&1&0&0&0 \\ 0&0&1&0&0 \\ 0&0&1&0&0  \end{array}\right)&=\left(\begin{array}{ccccc} 2&1&0&0&0\\1&0&1&0&0\\1&0&1&0&0\\0&1&1&0&0\\0&1&1&0&0\end{array}\right).
\intertext{Therefore the matrix of $[X,\cdot]: \nl_3^\sigma\rightarrow \nl_4^\sigma$ is}
\left(\begin{array}{ccc}2&1&0\\1&0&1\\0&1&1 \end{array}\right).&
\end{align*}
The determinant of this matrix is $-3=0$ and has at least rank $2$. Thus the cokernel of $[X,\cdot] : \nl_3^\sigma\rightarrow \nl_4^\sigma$ is $1$ dimensional.

\end{document}